\documentclass[12pt,a4paper]{article}

\usepackage{amsmath,amssymb,latexsym,amsthm}
\usepackage[english]{babel}
\usepackage[latin2]{inputenc}

\begin{document}

\newtheorem{Dfn}{Definition}
\newtheorem{Theo}{Theorem}
\newtheorem{Lemma}[Theo]{Lemma}
\newtheorem{Prop}[Theo]{Proposition}
\newtheorem{Coro}[Theo]{Corollary}
\newcommand{\Pro}{\noindent{\em Proof. }}
\newcommand{\Rem}{\noindent{\em Remark. }}

\title{$3$-dimensional loops on non-solvable reductive spaces}
\author{\'Agota Figula}
\date{}
\maketitle

\begin{abstract}  We treat the  almost differentiable left A-loops as images of global differentiable 
sharply transitive sections $\sigma :G/H \to G$ for a Lie group $G$ such that $G/H$ is a reductive homogeneous manifold. In this 
paper we classify all $3$-dimensional connected strongly left alternative almost differentiable left A-loops $L$, such that 
for the corresponding section $\sigma :G/H \to G$  the Lie group $G$ is non-solvable.    
\end{abstract}

\bigskip
\centerline{\bf Introduction}

\bigskip
\noindent
The associative law forces the identity  $(a b)^{-1} a b=1$ for all elements $a$  and $b$ of a group  $G$.  
For loops which are 
structures with a binary multiplication having up to associativity the same properties as groups this behaviour changes 
radically. This observation led to a broader research of loops $L$  in which the mapping 
$ x \mapsto  [(a b)^{-1}  (a (b x))]$ is an automorphism of $L$ (cf. \cite{bruck}, \cite{belousov}).  These loops have been called left A-loops. 

According to \cite{loops} we treat the left A-loops $L$ as images of global differentiable sharply transitive sections $\sigma :G/H \to G$ for a Lie group $G$ 
such that the subset $\sigma (G/H)$ is invariant under the conjugation with the elements of $H$. Here $G$  denotes  the group topologically generated by the left translations  of $L$ and $H$ is the stabilizer of the identity of $L$ in $G$. Loops given by a 
differentiable section in a Lie group are called almost differentiable. 

For an almost differentiable left A-loop $L$ 
the tangent space $T_1 \sigma (G/H)$ of the image of  $\sigma $ at $1 \in G$ can be provided with a binary and a ternary multiplication and  
yields a Lie triple algebra (cf. \cite{kikkawa1}, Definition 7.1, p. 173). Since the Lie triple algebras correspond to affine 
reductive spaces, which are  essential objects in  differential geometry (cf. \cite{foundations1}, 
\cite{helgason}), there is 
a strong connection between the theory of differential left A-loops  
 and the theory of affine reductive homogeneous spaces (cf.  \cite{kikkawa2}).
In particular the theory of  connected differentiable Bruck loops 
(which form a subclass of the class of left A-loops) is essentially the theory of affine 
symmetric spaces (cf. \cite{loops}, Section 11).      

The smallest dimension for a connected almost differentiable non-associa\-ti\-ve left A-loop is equal $2$. 
There exist precisely two isotopism classes of
$2$-dimensional  left A-loops. In the one class there lies only the hyperbolic plane loop which 
is related to the hyperbolic symmetric plane (cf. \cite{loops}, Section 22). 
In  the other isotopism class we may choose as a representative  the 
$2$-dimensional Bruck loop $L$ which  
is realized on the 
pseudo-euclidean affine plane $E$ such that the group topologically generated by its left translations  is the connected component of the group 
of pseudo-euclidean 
motions and the elements of $L$ are the lines of positive slope in E (cf. \cite{loops}, Section 25). 

Our aim in this paper is to classify the $3$-dimensional connected 
almost differentiable left A-loops, which have a non-solvable Lie group $G$ as the group topologically 
generated by their left translations and which correspond to differentiable sections $\sigma : G/H \to G$ such that the 
exponential image of the tangent space ${\bf m}=T_1(\sigma (G/H))$ is contained in $\sigma (G/H)$. These loops are called 
strongly left alternative almost differentiable left A-loops.  

Using the standard enveloping Lie algebra of a Lie triple algebra one sees that  $G$ is four, five or six dimensional. 
For the classification  we determine all complements ${\bf m}$ of the Lie algebra ${\bf h}$ of $H$ in the Lie algebra ${\bf g}$ 
of $G$ such that ${\bf m}$ generates ${\bf g}$ and satisfies the relation $[ {\bf h}, {\bf m}] \subseteq {\bf m}$. 
The submanifold $\exp {\bf m}$ can be extended to a global section if and only if $\exp {\bf m}$ forms a system of representatives for the cosets $\{ x H^g |\ x \in G \}$ in $G$, where $H^g=g^{-1} H g$ with $g \in G$.

In contrast to a frequent occurance of reductive spaces and hence strongly left alternative almost differentiable  local left A-loop, which 
can be 
represented as local sections in non-solvable Lie groups, the global loops in this class are rare and 
in strong relation to geometries on 3-dimensional manifolds as the following theorem shows.

\bigskip
{\bf Theorem} \begin{em} 
There are precisely two  
classes ${\cal C}_i\ (i=1,2)$ of  
connected almost differentiable strongly left alternative simple left A-loops $L$ having dimension $3$ such that the group $G$ 
generated by the left translations  is a non-solvable Lie group.

The  class ${\cal C}_1$ consists of left A-loops having the simple 
Lie group $G=PSL_2(\mathbb C)$ as the group topologically generated by their 
left translations, and the stabilizer $H$ of 
$e \in L$ in $G$ is the group $SO_3(\mathbb R)$. 
\newline
Any loop in the class ${\cal C}_1$ can be represented by a real 
parameter $a \in \mathbb R$.  For all $a \in \mathbb R$ the loops $L_a$ 
and $L_{-a}$ are isomorphic. These two loops form a full isotopism class.
The loops $L_a$, $a \in \mathbb R$ are realized on the hyperbolic symmetric space $H_3$ such that the group 
topologically generated by their left translations is the connected component of the group of motions of $H_3$. The 
elements of all loops $L_a$ in ${\cal C}_1$ are the points of $H_3$, but the sets of left translations differ.   
The hyperbolic 
space loop $L_0$, which is the unique Bruck loop in ${\cal C}_1$, is defined by the multiplication $x \cdot y= \tau_{e,x} (y)$, 
where $\tau_{e,x}$ is the hyperbolic translation 
moving $e$ onto $x$. 
\newline
The  class ${\cal C}_2$ of simple left A-loops  consists of  
$3$-dimensional connected differentiable left A-loops  
 such that the group  $G=PSL_2(\mathbb R) \ltimes \mathbb R^3$,  where 
the action of $PSL_2(\mathbb R)$ on $\mathbb R^3$ is the adjoint action of 
$PSL_2(\mathbb R)$ on its Lie algebra, is the group topologically generated 
by 
the left translations. This group is  the connected component of the group of pseudo-euclidean 
motions and the stabilizer $H$ of $e \in L$ in $G$ is the stabilizer of a plane on which the euclidean metric is 
induced.   
\newline
The loops in ${\cal C}_2$ can be represented by two real parameters $a,b$ and 
form precisely two isomorphism classes, which coincide with the isotopism 
classes. The  one isomorphism class consists of  Bruck loops $L_{a,0}$, $a \in \mathbb R$, and we may choose the 
pseudo-euclidean space loop $L_{0,0}={\hat L}_0$  as a representative of this isomorphism class. 
As a representative of the other 
isomorphism class which contains the loops $L_{a,b}$ with $b \neq 0$ may be chosen  
 the loop $L_{0,1}={\hat L}_1$. Any loop in the class  ${\cal C}_2$ is realized on the 
pseudo-euclidean affine 
space $E(2,1)$.  The elements of these loops are the planes on which the euclidean metric is induced but the sets of left translations differ.   
\newline
Moreover, the $3$-dimensional strongly left alternative almost differentiable  non-simple left A-loops 
are either the  products 
of a $1$-dimensional Lie group with a $2$-dimensional left A-loop isomorphic  to 
the hyperbolic plane loop or the  Scheerer extensions  of the Lie group 
$SO_2(\mathbb R)$ by the $2$-dimensional left A-loop isomorphic to the hyperbolic plane loop and the coverings of these Scheerer extensions.  
\end{em}

\bigskip
\noindent
Another class of almost differentiable loops which has been thoroughly 
investigated is the class of differentiable Bol loops. The sections 
$\sigma:G/H \to G$ of Bol loops are characterized by the fact that for all 
$a,b \in \sigma(G/H)$ the element $aba$ is also contained in $\sigma(G/H)$. 
The $3$-dimensional almost differentiable Bol loops with non-solvable Lie 
groups have been classified in \cite{figula}; the Lie groups $G$ 
topologically generated by their left translations as well as the 
corresponding stabilizers $H$ are the same as in the case of $3$-dimensional 
almost differentiable left A-loops, but the sections  essentially differ. The 
intersection of these two classes are only the Bruck loops and the Scheerer 
extensions of the orthogonal group $SO_2(\mathbb R)$ by the hyperbolic plane loop and the coverings of these Scheerer 
extensions.

\bigskip
\centerline{\bf 1. Left A-loops }

\bigskip
\noindent
{\bf 1.1} A set $L$ with a binary operation $(x,y) \mapsto x \cdot y$ is called a loop 
if there exists an element $e \in L$ such that $x=e \cdot x=x \cdot e$ holds 
for all $x \in L$ and the equations $a \cdot y=b$ and $x \cdot a=b$ have 
precisely one solution which we denote by $y=a \backslash b$ and $x=b/a$. 
The left 
translation $\lambda _a: y \mapsto a \cdot y :L \to L$ is a bijection of $L$ 
for any $a \in L$.
Two loops $(L_1, \circ )$ and $(L_2, \ast )$ are called isotopic if there are 
three bijections $\alpha ,\beta ,\gamma : L_1 \to L_2 $ such that 
$\alpha (x) \ast \beta (y)=\gamma (x \circ y)$ holds for any $x,y \in L_1$. 
An isotopism is an equivalence relation. 
If $\alpha =\beta =\gamma $ then the isotopic loops $(L_1, \circ )$ and 
$(L_2, \ast )$ are called isomorphic.
Let $(L_1, \cdot )$ and 
$(L_2, \ast )$ be two loops. The direct product $L=L_1 \times L_2= \{ (a,b)
\ |a \in L_1, b \in L_2 \}$ with the multiplication 
$(a_1,b_1) \circ (a_2,b_2)=(a_1 \cdot a_2, b_1 \ast b_2)$ is again a loop, 
which is called the direct product of $L_1$ and $L_2$, and the loops 
$(L_1, \cdot )$,  $(L_2, \ast )$ are subloops of $(L, \circ )$.   
\newline
A loop is called a left A-loop if each mapping 
$\lambda _{x,y}=\lambda _{xy}^{-1} \lambda _x \lambda _y:L \to L$ is an 
automorphism of $L$.
\newline
Let $G$ be the group generated by the left 
translations of $L$ and let $H$ be the stabilizer of $e \in L$ in the group 
$G$. 
The left translations of $L$ form a subset of $G$ acting on the cosets 
$\{x H; x \in G\}$ such that for any given cosets $aH$ and $bH$ there exists 
precisely one left translation $ \lambda _z$ with $ \lambda _z a H=b H$. 
\newline
Conversely let $G$ be a group, H be a subgroup containing no normal 
non-trivial subgroup of $G$ and $\sigma : G/H \to G$ be a section with $\sigma (H)=1 \in G$ such that the set $\sigma (G/H)$ of 
representatives for the left cosets $\{ x H, x \in G \}$  
 and acts sharply transitively on the space $G/H$ of  $\{x H, x \in G\}$ (cf. \cite{loops}, p. 18). Such a section we call a sharply transitive section. 
Then the multiplication  defined by 
$x H \ast y H=\sigma (x H) y H$ on the factor space $G/H$ or by $x \ast y=\sigma(xyH)$ on $\sigma (G/H)$ yields a loop $L(\sigma )$. The group $G$ is isomorphic to the group generated by the left translations of $L(\sigma )$. 
\newline
If $G$ is a Lie group and $\sigma $ is a differentiable section 
satisfying the above conditions then the loop $L(\sigma )$ is almost differentiable. 
This loop is a left A-loop if and only if the subset 
$\sigma (G/H)$ is invariant under the conjugation with the elements of $H$.
\newline
Let $L_1$ be a loop defined on the factor space $G_1/H_1$ with respect to a 
section 
$\sigma _1:G_1/H_1 \to G_1$ the image of which is the set $M_1 \subset G_1$. 
Let $G_2$ be a  group let $\varphi :H_1 \to G_2$ be a homomorphism and $(H_1, \varphi (H_1))=\{(x, \varphi (x)); x \in H_1\}$. A 
loop $L$ is called a Scheerer extension of 
$G_2$ by $L_1$ if the loop $L$ is defined on the factor space 
 $(G_1 \times G_2)/(H_1, \varphi (H_1))$   
with respect to the section $\sigma:(G_1 \times G_2)/(H_1, \varphi (H_1)) \to G_1 \times G_2$ the image of which is the set $M_1 \times G_2$. 
\newline
The loops $L_1$ and $L_2$ having the same group $G$ of the group generated by the left translations 
and the same stabilizer $H$ of $e \in L_1, L_2$ are isomorphic if there is 
an automorphism of $G$ leaving $H$ invariant and mapping the section 
$\sigma_1(G/H)$ onto the section $\sigma_2(G/H)$. Moreover let $L$ and $L'$ 
be loops having the same group $G$ generated by their left translations. 
Then $L$ and $L'$ are isotopic if and only if there is a loop $L''$ 
isomorphic to $L'$ having $G$ again as the group generated by its left 
translations such that there exists an inner automorphism $\tau $ of $G$ 
mapping the section $\sigma ''(G/H)$ belonging to $L''$ onto the section 
$\sigma (G/H)$ corresponding to $L$ (cf. 
\cite{loops}, Theorem 1.11.  pp. 21-22). 
\newline
If $L$ is a connected almost differentiable left A-loop, then the group $G$ 
topologically generated by the left translations of $L$ within the group of 
autohomeomorphisms is
 a connected Lie group (cf. \cite{quasigroups}; \cite{loops}, Proposition 5.20. p. 75), and we may describe
 $L$ by a 
differentiable section. 
\newline
Let $L$ be a connected almost differentiable left A-loop. Let $G$ be the Lie 
group topologically generated by the left translations of $L$, and let 
$(\bf{g, [.,.]}) $ be the Lie algebra of 
$G$. Denote by $\bf{h}$ the Lie algebra of the stabilizer $H$ of the identity $e \in L$ 
in 
$G$ and by ${\bf m}=T_1 \sigma (G/H)$ the tangent space at $1 \in G$ of 
the image of the section $\sigma :G/H \to G$ corresponding to 
 $L$. Then ${\bf m}$ generates ${\bf g}$ and the homogeneous space $G/H$ is 
reductive, i.e.
we have $\bf{g}=\bf{m} \oplus \bf{h}$ and 
$[{\bf h}, {\bf m}] \subseteq {\bf m}$. 
(cf. \cite{foundations1} Vol II, p. 190; \cite{loops}, Proposition 5.20. p. 75)  
The subspace ${\bf m}$  with the operations $X \cdot Y:=[X,Y]_{{\bf m}}$ and $[X,Y,Z]:=[[X,Y]_{{\bf h}},Z ]$ 
yields a Lie triple algebra (\cite{kikkawa1}, Definition 7.1, p. 173).  
If $X \cdot Y=0$ for all $X,Y \in {\bf m}$ then ${\bf m}$ is a Lie triple system.
In this case  the factor space $G/H$ is an affine symmetric space (\cite{symmetric})  and the corresponding loop $L$ is 
called a Bruck loop. The Lie algebra ${\bf g}$ of $G$ is isomorphic to the standard enveloping Lie algebra of the Lie triple 
algebra ${\bf m}$ generating ${\bf g}$. 
If the dimension of ${\bf m}$ is  $n$ then  ${\bf g}$ has dimension  at most $n(n+1)/2$. 
\newline
In this paper we investigate  strongly left alternative (cf. \cite{loops}, Definition  5.3, p. 67) almost 
differentiable  left A-loops $L$ of  dimension $3$; these loops satisfy  
 $\exp [T_1 \sigma (L)] \subset \sigma (L)$ (\cite{loops}, Proposition 5.5 p. 68).  
Hence every global left 
A-loop contains an 
exponential image of a 
complement ${\bf m}$ of the Lie algebra ${\bf h}$ of $H$ in the Lie algebra 
${\bf g}$ of $G$, such that ${\bf m}$ generates ${\bf g}$ and satisfies the 
relation $ [{\bf h},{\bf m}] \subseteq {\bf m}$. 
\newline
In this paper we often compute the images of  subspaces ${\bf m}$ of the Lie algebras $sl_2(\mathbb R)$, 
$sl_2(\mathbb C)$, $su_2(\mathbb C)$ under the exponential map. 
\newline
{\bf 1.2} The exponential function of the Lie algebras $sl_2(\mathbb R)$, 
$sl_2(\mathbb C)$, $su_2(\mathbb C)$.  
\newline
The exponential map $\exp : {\bf g} \to G$ is defined in the following way: For
$X \in {\bf g}$ we have $\exp \ X= \gamma _X(1)$, where $\gamma _X(t)$ is the 
1-parameter subgroup of $G$ with the property 
$ \frac{d}{dt} \big |_{t=0} \gamma _X(t)=X$. 
The matrices 
\newline
\centerline{$ K= \left ( \begin{array}{rr}
1 & 0 \\
0 & -1 \end{array} \right )$, \ \ $T=\left ( \begin{array}{cc}
0 & 1 \\
1 & 0 \end{array} \right )$, \ \ $U=\left ( \begin{array}{cc}
0 & 1 \\
-1 & 0 \end{array} \right )$ } 
form  a real basis of $sl_2(\mathbb R)$.  
The Lie algebra multiplication is given by the rules:
\newline
\centerline{$ [K,T]=2 U, \ \  [K,U]=2 T, \ \ [U,T]=2 K$ . }
The normalized Cartan-Killing form $k:sl_2(\mathbb R) 
\times sl_2(\mathbb R) \to \mathbb R$ of $sl_2(\mathbb R)$  
 is the 
bilinear form
defined by  $k(X,Y)=\displaystyle \frac {1}{8} 
\hbox{trace} (\hbox{ad}X\ \hbox{ad}Y)$. 
\newline
If $X \in sl_2(\mathbb R)$ 
has the decomposition 
\newline
\centerline{
$X=\lambda _1\  K+\lambda_2\ T+
\lambda_3\  U$} then the  Cartan-Killing form
$k$ satisfies 
\newline
\centerline{
$k(X)=\lambda _1^2+\lambda_2^2- \lambda_3^2 $.}
According to $\cite{old}$  for the exponential function  $\exp :sl_2(\mathbb R) \to SL_2(\mathbb R)$ we have 
\[ \exp\ X=C(k(X))\ I+S(k(X))\ X. \]   
Here is 
\newline
\centerline{$C(x)=  \begin{array}{c}
\cosh \sqrt{x} \quad \hbox{for} \ \ \ 0 \le x, \\
\cos \sqrt{-x} \quad \hbox{for} \ \ \ 0 > x, \end{array} \quad 
\sqrt{ | x |} \ S(x)=  \begin{array}{c}
\sinh \sqrt{x} \quad \hbox{for} \ \ \ 0 \le x, \\
\sin \sqrt{-x} \quad \hbox{for} \ \ \ 0 > x. \end{array}  $}
As a natural generalization of this formula 
we obtain the explicite form for the exponential function of 
$sl_2(\mathbb C)$.
Representing the Lie algebra ${\bf g}=sl_2(\mathbb C)$ as complex $(2 \times 2)$-matrices we may choose as basis $\{ K, T, U, \hbox{i} K, \hbox{i} T, \hbox{i} U \}$, where  $K, T, U$ are the basis elements of $sl_2(\mathbb R)$ (see in {\bf 1.2}).  
\newline
The normalized complex Cartan-Killing form $k_{\mathbb C}:sl_2(\mathbb C) 
\times sl_2(\mathbb C) \to \mathbb C$ of $sl_2(\mathbb C)$  
 is the 
bilinear form
defined by: $k_{\mathbb C}(X,Y)=\displaystyle \frac {1}{8} 
\hbox{trace} (\hbox{ad}X\ \hbox{ad}Y)$. If $X \in sl_2(\mathbb C)$ 
has the decomposition 
\newline
\centerline{
$X=\lambda _1\  K+\lambda_2\  T+
\lambda_3\  U+ \lambda_4\  \hbox{i} K+ \lambda_5\  \hbox{i} 
 T +\lambda_6\  \hbox{i}  U$} then the complex Cartan-Killing form
$k_{\mathbb C}$ satisfies 
\newline
\centerline{
$k_{\mathbb C}(X)=\lambda _1^2+\lambda_2^2+
\lambda_6^2-
\lambda_3^2- \lambda_4^2- \lambda_5^2+i\ (2\ \lambda _1  \lambda_4+2\ 
\lambda_2 \lambda_5-2\ \lambda_3 \lambda_6)$ } (cf. \cite{freudenthal}, Section 1, pp. 1-3).
The normalized real Cartan-Killing form $k_{\mathbb R}:sl_2(\mathbb C) \times 
sl_2(\mathbb C) \to \mathbb R$ is the restriction of $k_{\mathbb C}$ to 
$\mathbb R$ such that 
\newline
\centerline{
$k_{\mathbb R}(X)=\lambda _1^2+\lambda_2^2+\lambda_6^2-
\lambda_3^2- \lambda_4^2- \lambda_5^2. $} 
\newline
For the exponential function $\exp:sl_2(\mathbb C) \to SL_2(\mathbb C)$ one has 
\[ \exp \  {X}=C(k_{\mathbb C}(X))\ I+S(k_{\mathbb C}(X))\ X, \]
where $C(z)=\cosh \sqrt{z}$ and $S(z)=\displaystyle \frac{\sinh \sqrt{z}}{\sqrt{z}}$, $z \in \mathbb C$.
\newline
The group $SU_2(\mathbb C)$ is the $3$-dimensional compact subgroup of $SL_2(\mathbb C)$, which can be represented by $(2 \times 2)$-complex matrices 
of the form:
\newline
\centerline{
$\left \{ \left ( \begin{array}{rr}
a & b \\ 
-\bar{b} & \bar{a}
\end{array} \right); a,b \in \mathbb C, a\bar{a}+b\bar{b}=1 \right \}. $} 
Therefore the Lie algebra ${\bf g}=su_2(\mathbb C)$ is generated by the basis elements 
$U$, $i K$, $i T$. 
The restriction of the  formula for the exponential function of  
$sl_2(\mathbb C)$ to  $su_2(\mathbb C)$ gives the  formula for $\exp : su_2(\mathbb C) \to SU_2(\mathbb C)$. 
\newline
\newline
{\bf 1.3} Now we study which pairs $(G,H)$ of Lie groups can  admit differentiable sections $\sigma :G/H \to G$ 
corresponding to  $3$-dimensional almost differentiable left A-loops.
\newline
We start with a well 
known fact from linear algebra: 
\begin{Lemma}
If $\bf{g}= \bf{a} \oplus \bf{b}$ with a $3$-dimensional subspace ${\bf a}$ and the dimension of ${\bf g}$ is  4 or 5 then  
 ${\bf m} \cap {\bf a}$ is at least $2$ respectively $1$-dimensional for any $3$-dimensional subspace 
${\bf m}$. 
\end{Lemma}

\noindent
The next fact is proved in \cite{figula}, Lemma 3. 
\begin{Lemma} Let $L$ be an almost differentiable global loop and denote by 
${\bf m}$ the 
tangent space of $T_1 \sigma (G/H)$, where $\sigma :G/H \to G$ is the section corresponding to $L$. Then ${\bf m}$ does not 
contain any 
 element of $Ad_{g^{-1}} {\bf h}=g {\bf h} g^{-1} $ for some $g \in G$. Moreover every element of $G$ can be uniquely written  as 
a product of an element of $\sigma (G/H)$ with an element of $H$.  
\end{Lemma}

\noindent
Since a $1$-dimensional almost differentiable left A-loop is a group an analogue of Proposition 1 in 
\cite{figula} is the following 
\begin{Prop}
Let $L$ be a loop
and let $G$ be the group generated by the left translations of $L$, and 
denote by
$H$ the stabilizer of $e \in L$ in $G$. If $G$ and $H$ are direct products 
$G=G_1 \times G_2$ and $H=H_1 \times H_2$ with $H_i \subset G_i$ $(i=1,2)$ 
then  $L$ is the  product 
of two loops $L_1$ and $L_2$, 
and $L_i$ is isomorphic to a loop $L_i^{ \ast }$ having $G_i$ as the group 
generated by the left translations of $L_i^{ \ast }$ and $H_i$ as the 
corresponding stabilizer subgroup $(i=1,2)$.  
\newline
In particular there exists no $3$-dimensional left A-loop $L$ such that $L$
is the  product of a $1$-dimensional and a $2$-dimensional left A-loop 
and $L$ has a $5$- or $6$-dimensional Lie group as the group topologically 
generated by its left translations. 
\end{Prop}

\begin{Lemma} Let ${\bf g}={\bf g}_1 \oplus {\bf g}_2$ be the Lie algebra of the Lie group 
$G=G_1 \times G_2$, 
such that $G$ is the group topologically generated by the left translations of a $3$-dimensional almost 
differentiable left A-loop $L$. Let ${\bf m}$ be the tangent space of the manifold $\Lambda $ of the left 
translations of 
$L$ at $1 \in G$. Denote by ${\bf h}$ the Lie algebra of the stabilizer $H$ of $e \in L$ in $G$ and let 
$\pi _i: {\bf g} \to {\bf g}_i$, $i=1,2$ be the natural projection of ${\bf g}$ onto ${\bf g}_i$. We assume that
 ${\bf g}_1$ is isomorphic to $sl_2(\mathbb R)$  and $\hbox{dim}\ \pi_1 ({\bf h})=2$. Then: 
\newline
(i) $\hbox{dim}\ \pi_1 ({\bf m})=3$. 
\newline
(ii) If $\hbox{dim}\ \pi_2 ({\bf m}) \ge 2$ then the Lie algebra 
${\bf h}$ has the form 
\newline
\centerline{ ${\bf h}=\{ (x, \varphi (x))| x \in \pi_1({\bf h}) \}$, }
with an isomorphism $\varphi :\pi_1({\bf h}) \to \pi_2({\bf h})$. Moreover, one has  
$\hbox{dim}\ \pi_2 ({\bf m})= 2$ and $\hbox{dim}\ {\bf g}=5$. 
\end{Lemma}
\Pro (i) We have      
$\hbox{dim}\ \pi_1 ({\bf m}) \ge 2$ since otherwise the set $\Lambda $ would not generate $G$. If  
$\hbox{dim}\ \pi_1 ({\bf h})=2$ then we may assume that  $\pi _1({\bf h})$ is generated by the elements $K, U+T$ 
of the Lie 
algebra $sl_2(\mathbb R)$ (see {\bf 1.2} in section 1). If $\pi_1({\bf m})$ were $2$-dimensional then 
it has one of the following forms (up to conjugation)  
\newline
a) $\pi_1({\bf m})= \langle U+a_1 K +a_2 (U+T), K+ b_1 K+b_2 (U+T) \rangle$,  
\newline
b) $\pi_1({\bf m})= \langle U+a_1 K +a_2 (U+T), T+ b_1 K+b_2 (U+T) \rangle$, 
\newline
c) $\pi_1({\bf m})= \langle T+a_1 K +a_2 (U+T), K+ b_1 K+b_2 (U+T) \rangle$, 
where $a_1,a_2,b_1,b_2 \in \mathbb R$. 
One has 
\newline
($\ast $) \centerline{$ \pi_1([{\bf h}, {\bf m}])=[\pi_1({\bf h}), \pi_1({\bf m})] \subseteq \pi_1({\bf m})$. }
In the case a) the element 
\newline
\centerline{$[K, U+a_1 K +a_2 (U+T)]=2 T+2 a_2(U+T) $}
is contained in $\pi _1({\bf m})$ if and only if $a_1=0, a_2=- \frac{1}{2}$. 
Moreover the element $[U+T, U-T]=-4 K$ is contained in $\pi _1({\bf m})$ precisely if $b_2 =0$ and $b_1 \neq -1$. 
But then $\pi_1({\bf m})=\langle (1+b_1)K, U-T \rangle $ is a subalgebra of ${\bf g}_1$.
\newline
In the case b) the element 
\newline
\centerline{$[U+T, U+a_1 K +a_2 (U+T)]=2 K+2 a_1(U+T) $} is not contained in  $\pi _1({\bf m})$, this is a 
contradiction to ($\ast $).
\newline
In the case c) we obtain the same contradiction in the same way as in the case a). Therefore is 
$\hbox{dim}\ \pi _1({\bf m})=3$. 
\newline
(ii) If $\hbox{dim}\ \pi _2({\bf h})=3$ then one has $ \pi _2({\bf h})={\bf g}_2$ and ${\bf h} \cap (0,{\bf g}_2) 
\neq (0,0)$. Then there exists a homomorphism $\beta : \pi _2({\bf h}) \to  \pi _1({\bf h})$ such that 
${\bf h} \cap  (0,{\bf g}_2)=\beta ^{-1}(0)$. This is a contradiction since ${\bf h}$ does not contain non-trivial ideal of ${\bf g}$.
\newline
Since ${\bf m} \subseteq  \pi _1({\bf m}) \times  \pi _2({\bf m})$ and according (i)  
$\hbox{dim}\ \pi _1({\bf m})=3$ there is a linear mapping $\alpha :  \pi _1({\bf m}) \to  \pi _2({\bf m})$ such that $\alpha $ is a linear isomorphism if $\hbox{dim}\ \pi _2({\bf m})=3$ and $\hbox{dim}\ \alpha ^{-1}(0)=1$ for 
$\hbox{dim}\ \pi _2({\bf m})=2$. 
\newline 
If $\hbox{dim}\ \pi _2({\bf h}) \le 2$ and (ii) does not hold then there is a homomorphism 
$\gamma : \pi _1({\bf h}) \to  \pi _2({\bf h})$ with $0 \neq S=Ker\ \gamma $. If  $\hbox{dim}\  S=2$ then we have 
$\pi _1({\bf h})=S$ and the Proposition $3$ gives a contradiction. Hence $\hbox{dim}\  S=1$. Then 
$(S,0)={\bf h} \cap (\pi_1({\bf m}), 0)={\bf h} \cap ({\bf g}_1, 0)$ and $(S,0)=\langle (U+T, 0) \rangle $ is a 
$1$-dimensional subalgebra of ${\bf h}$. 
\newline
First we treat the case that $\hbox{dim}\  \pi _2({\bf m})= 3$. Then one has ${\bf m} \cap ({\bf g}_1,0)=(0,0)$ and 
$\pi_2({\bf m}_2)={\bf g}_2$. Then there is an element $m_2 \in \pi_2({\bf m})$ such that 
\newline
\centerline{$[(r\ (U+T), 0), (\alpha (m_2), m_2)]=([r\ (U+T),\alpha (m_2)],0) \neq (0,0)$, }
where $r \in \mathbb R$. This is a contradiction.  
\newline
Now we assume $\hbox{dim}\  \pi _2({\bf m})=2$. Then we have ${\bf m} \cap ({\bf g}_1,0)=(S',0)$, where $S'=\alpha ^{-1}(0)$.
Since 
\newline
\centerline{$[{\bf h}, {\bf m} \cap ({\bf g}_1,0)]=[(h_1,h_2), (m_1,0)]=([h_1,m_1],0) \subset  {\bf m} \cap ({\bf g}_1,0)$ }
with $(h_1,h_2) \in {\bf h}$ and $(m_1,0) \in {\bf m}$ it follows that $\pi_1({\bf h})$ normalizes $S'$ and 
therefore $(S',0)=\langle (U+T,0) \rangle =(S,0)$, which is a contradiction. 
\qed

\begin{Prop} Let $G=G_1 \times G_2$ be the group 
topologically generated by the left translations of a $3$-dimensional 
connected almost differentiable proper left A-loop $L$.  Let the group $G_1$ be locally isomorphic either 
to $SO_3(\mathbb R)$ or to $PSL_2(\mathbb R)$. Then for the pair $(G,H)$, where $H$ is the stabilizer 
of $e \in L$ in $G$, one of the following cases occurs: 
\newline
1)  $L$ is the   
product of the hyperbolic plane loop with a $1$-dimensional Lie group  and 
$H \cong SO_2(\mathbb R) \times \{ 1 \}$. 
\newline
2) $G$ is isomorphic to $PSL_2(\mathbb R) \times \mathbb R$ and 
$H=\{ ( x, \varphi (x))\}$, where $\varphi $ is a monomorphism from the 
$1$-dimensional subgroup $\left \{ \left ( \begin{array}{cc} 
1 & b \\
0 & 1 \end{array} \right ), b \in \mathbb R \right \}$ or from 
$\left \{ \left ( \begin{array}{ll}
a & 0 \\
0 & a^{-1} \end{array} \right ), a > 0 \right \}$ of $PSL_2(\mathbb R)$ 
onto $\mathbb R$. 
\newline
3) $G=PSL_2(\mathbb R) \times SO_2(\mathbb R)$ such that 
$H=\{ ( x, x^n) | x \in SO_2(\mathbb R), n \in \mathbb N \}$.  
\newline
4) $G \cong  PSL_2(\mathbb R) \times PSL_2(\mathbb R)$ and $H$ has the form 
\newline
\centerline{$H=\{ ( x,x )\ |\  x \in PSL_2(\mathbb R) \}$. } 
\end{Prop}
\Pro First we assume that $\hbox{dim}\ G=6$. If $G_1$ is locally isomorphic to $PSL_2(\mathbb R)$ 
then it 
follows from Lemma 4 and from the proof of Proposition 4 in \cite{figula} that we are in the case 4). If 
$G_1$ is locally isomorphic to  $SO_3(\mathbb R)$ then it is easy to see that  $G_2$ is also locally isomorphic 
to  $SO_3(\mathbb R)$ and we may assume that 
$H=\{ ( x,x )\ |\  x \in G_1 \}$. This case is excluded by 
Proposition 16.11 in \cite{loops} (p. 205). 
\newline
Now we suppose that  $\hbox{dim}\ G=5$. We may assume that $\hbox{dim}\  \pi _1({\bf h})= 2$ since 
otherwise $H$ would be a direct product $H=H_1 \times H_2$ which contradicts Proposition 3. Now it 
follows from Lemma 4 that 
\newline
\centerline{$H= \{ (x, \varphi (x))| x \in \pi_1(H) \}$,} 
where $\pi_1 (H)$ is isomorphic to the group 
${\cal L}_2=\{ x \mapsto a x + b; a >0, b \in \mathbb R \}$ and  $\varphi :\pi_1(H) \to \pi_2(H) $ 
is an isomorphism. A real basis of the Lie algebra
$\bf{g}=sl_2(\mathbb R) \oplus {\cal L}_2$ is  
\newline
\centerline{${\bf g}=\langle (K,0),\ (T,0),\ (U,0),\ (0,e_1),\ 
(0,e_2) \rangle $,} where $K,T$ and $U$ are the basis elements of $sl_2(\mathbb R)$ 
(see {\bf 1.2}) and $e_1,e_2$ are the basis elements of ${\cal L}_2$ with the rule 
$[(0, e_1),(0, e_2)]=-(0, e_2)$.  The Lie algebra ${\bf h}$ of $H$ is given by 
\newline
\centerline{${\bf h}=\langle (K,e_1), (U+T,e_2) \rangle $.}
An arbitrary complement ${\bf m}$ to ${\bf h}$  
 in ${\bf g}$ has as generators  
\newline
\centerline{$l_1=(U+a_1 K +a_2 (U+T), a_1 e_1 + a_2 e_2)$, } 
\newline
\centerline{$l_2=(b_1 K +b_2 (U+T), e_1 +b_1 e_1+b_2 e_2)$, }
\newline
\centerline{$l_3=(c_1 K +c_2 (U+T), e_2 + c_1 e_1 +c_2 e_2)$, }
where $a_1,a_2,b_1,b_2,c_1,c_2 \in \mathbb R$. 
Since one has $\hbox{dim} \ {\bf m} \cap (sl_2(\mathbb R) \oplus \{ 0\}) \ge 1$ and 
the relation $[ {\bf h}, {\bf m} ] \subseteq {\bf m} $ holds we obtain that $G$ cannot have dimension 5. 
\newline
Finally let $\hbox{dim}\ G=4$. Then $\hbox{dim}\ H=1$. First we assume that $G_1$ is locally isomorphic 
to  $PSL_2(\mathbb R)$. If $\pi_2(H)=1$ according to Proposition 3 and Theorem 27.1, Theorem 18.14 in 
\cite{loops} one has $H=SO_2(\mathbb R) \times \{ 1 \}$ and $L$ is a product of 
 the hyperbolic plane loop with a $1$-dimensional Lie group. This is the case 1).  
\newline
Let now   $\pi_2(H) \neq 1$.  If  $G_2$ is isomorphic to $\mathbb R$ then $H=\{ (\varphi (x), x)\ |\ x \in 
\mathbb R \}$, where  $\varphi $ is a  monomorphism onto $G_1$. The inverse of $\varphi $ is again a 
monomorphism. Since the group 
$PSL_2(\mathbb R)$ has precisely 2 conjugacy classes of $1$-dimensional subgroups isomorphic to 
$\mathbb R$  we obtain the cases  2). If $G_2$ is isomorphic to $SO_2(\mathbb R)$ then we may assume 
that 
\newline
\centerline{$H=\{ (x, x^n) | x \in SO_2(\mathbb R), n \in \mathbb N  \}$. }  
It remains to consider a group $G$ locally isomorphic to $SO_3(\mathbb R) \times SO_2(\mathbb R)$. 
Since $H$ 
does not contain any non-trivial normal subgroup the group $G$ is isomorphic to   
$SO_3(\mathbb R) \times SO_2(\mathbb R)$ and $H$ has one of the following forms: 
\newline
\centerline{$H'=\{ K \times \{ 0 \} \}$, \ \ or \ \   $H=\{ ( k, \varphi (k)) | k \in K \}$, }
where $K$ is isomorphic to $SO_2(\mathbb R)$ and $\varphi $ is a non-trivial homomorphism. 
\newline
Since in the first case the factor space $G/H'$ is a topological product of spaces having as a factor 
the 
$2$-sphere or the projective plane we have to consider only the second case (\cite{loops}, 
Theorem 19.1, p. 249).
\newline   
The Lie algebra ${\bf g}$ of $G$ can be represented as $su_2(\mathbb C) \oplus \mathbb R$. 
Then as a basis of ${\bf g}$ may be chosen the following elements  
$i (K,0), (U,0), i (T,0)$, $(0,e_1)$,  
where $i K, U, i T$ is the real basis of $su_2(\mathbb C)$ which is introduced in {\bf 1.2} and $e_1$ is the basis element of $\mathbb R$. 
Moreover the Lie group $H$ has one of the following shapes  
$H_n= \{ (x, x^n) | x \in SO_2(\mathbb R), n \in \mathbb N  \}$ and 
for the Lie algebra ${\bf h}$ of $H_n$ has the form  ${\bf h}=\langle (U, e_1) \rangle$.  
An arbitrary complement ${\bf m}$ to the Lie algebra ${\bf h}$ of $H_n$ in ${\bf g}$ has the shape:
\[ {\bf m}= \langle (i K + a_1\ U,a_1\  e_1), (i T + a_2\  U, a_2\  e_1), 
(a_3\  U, e_1 + a_3\  e_1) \rangle ,\]
where $a_1, a_2, a_3 \in \mathbb R$.
From Lemma 1  and from the property $[{\bf h}, {\bf m}] \subseteq {\bf m}$ we obtain that  
the unique reductive complement ${\bf m}$ generating ${\bf g}$ has the form 
\newline
\centerline{ ${\bf m}=\langle i(K,0), i(T,0), (a\ U, e_1 + a\  e_1) \rangle , $}
where $a \in \mathbb R \backslash \{ -1 \}$. 
\newline
For $a >-\displaystyle \frac{1}{2}$ and the basis element $(U,e_1) \in {\bf h}$ one has 
\newline
\centerline{
$Ad _g (U,e_1)=-2 k l (i T,0)+ ( \frac{a}{1+a} U,e_1) \in {\bf m}_{a}$} 
with $g= \left ( \pm \left (\begin{array}{cc}
k- l i & 0  \\
0 & k+l i  \end{array} \right), 0  \right) \in G$, such that $k^2-l^2 = \frac{a}{1+a}$ and $k^2+l^2=1$.  This contradicts Lemma 2.  
\newline
For $a < - \frac{1}{2}$ the vectors 
\newline
\centerline{
 $X_1= \left ( \frac{\pi}{6} U + \frac{\sqrt{143}}{6} \pi i K, \frac{\pi (1+a)}{6 a} e_1 \right ) \quad \hbox{and} $}
\newline
\centerline{
$X_2= \left ( \frac{\pi (1+a)}{6(1+a-na)} U, \frac{\pi (1+a)^2}{6 a(1+a-na)} e_1 \right ) $} 
are contained in ${\bf m}_a$. According to {\bf 1.2} we get 
\[ \exp X_1= \left ( \pm I,  \left ( \begin{array}{rr}
\cos \frac{\pi (1+a)}{6 a} &  \sin \frac{\pi (1+a)}{6 a} \\
-\sin \frac{\pi (1+a)}{6 a} & \cos \frac{\pi (1+a)}{6 a} \end{array} 
\right ) \right), \] 
\[ \exp X_2=  \left ( \pm \left ( \begin{array}{rr}
\cos l &  \sin l \\
-\sin l & \cos l \end{array} 
\right ) ,  \left ( \begin{array}{rr}
\cos \frac{l (1+a)}{ a} &  \sin \frac{l (1+a)}{ a} \\
-\sin \frac{l (1+a)}{ a} & \cos \frac{l (1+a)}{ a} \end{array} 
\right ) \right), \]  where $l= \frac{\pi (1+a)}{6(1+a-na)}$, $\pm I$ is the identity of $SO_3(\mathbb R)$. 
For the element 
\newline
\centerline{
$g= \left ( \pm I , \left ( \begin{array}{rr}
\cos \frac{\pi (1+a)}{6 a} &  \sin \frac{\pi (1+a)}{6 a} \\
-\sin \frac{\pi (1+a)}{6 a} & \cos \frac{\pi (1+a)}{6 a} \end{array} 
\right ) \right) \in G$} 
one has 
\newline
\centerline{$ g= \exp X_1= \exp X_2 \cdot h $}
with 
\newline
\centerline{$ h= \left (\pm \left ( \begin{array}{rr}
\cos l &  -\sin l \\
\sin l & \cos l \end{array} 
\right ) , \left ( \begin{array}{rr}
\cos  n l & -\sin n l \\
\sin n l & \cos n l \end{array} \right ) \right). $} 
This is again a contradiction to  Lemma 2. Therefore there is no global section $\sigma :G/H_n \to G$ satisfying $\exp {\bf m}_a 
\subseteq \sigma(G/H_n)$. 
\qed

\begin{Coro}
There is no global left A-loop $L$   
homeomorphic to the compact space $S^3$ or $P^3$. 
\end{Coro}
\Pro  The group $G$ topologically generated by the left translations of an almost differentiable proper left A-loop $L$ homeomorphic to $S^3$ acts transitively on $L$. According to 96.16 in $\cite{salzmann2}$ any maximal compact subgroup of $G$ acts also transitively on $S^3$. Since a transitive compact 
subgroup of $G$ is a non-solvable subgroup of $SO_4(\mathbb R)$ (96.20 in 
$\cite{salzmann2}$) the group $G$ is non-solvable. According to Proposition 16.11 in $\cite{loops}$ and  Proposition 5  there is no almost differentiable left A-loop homeomorphic to $S^3$ or $P^3$ having a non-solvable Lie group as the group topologically generated by its left translations. 
\qed

\bigskip
\centerline{\bf 2. Left A-loops as sections in semisimple Lie groups}

\bigskip
\noindent
In this section we classify all 3-dimensional connected strongly left alternative almost differentiable left A-loops $L$ having 
a semisimple Lie group $G$ as the group topologically generated by its left translations  and  
describe the reductive spaces and natural geometries associated with them.  
\newline
It follows from Lemma 4 and Proposition 5 that the group $G$  must be locally isomorphic either to $PSL_2(\mathbb C)$ or to 
$PSL_2(\mathbb R) \times PSL_2(\mathbb R)$. In the second case we may assume that the stabilizer $H$ of $e \in L$ in $G$  is locally 
isomorphic to $\{ (x,x); x \in PSL_2(\mathbb R) \}$.  

\begin{Lemma} For all $\lambda \in \mathbb R \backslash \{ 0,1 \}$ there is a reductive complement 
\newline
\centerline{
${\bf m}_{\lambda }= \{ ( X, \lambda X) | X \in sl_2(\mathbb R ) \}$}
 to the Lie algebra ${\bf h}=\{ (X,X);\ X \in sl_2 (\mathbb R) \}$ of 
$H$ in ${\bf g}=sl_2(\mathbb R) \oplus sl_2(\mathbb R)$. 
\end{Lemma}
\Pro Let $K, U$ and $T$ be the real basis of $sl_2(\mathbb R)$ induced in {\bf 1.2}. 
In this case a $3$-dimensional  complement ${\bf m} \subset {\bf g}$ has the shape 
\newline
\centerline{
$\{(X, \varphi (X))| X \in sl_2(\mathbb R) \}$,}
 where $\varphi : sl_2(\mathbb R) \to sl_2(\mathbb R)$ is a linear map. 
From the relation $[{\bf h}, {\bf m}] \subseteq {\bf m}$  we obtain the assertion.  
\qed 

\begin{Prop} There is no global sharply transitive section $\sigma :G/H \to G$ satisfying the relation 
$\sigma (G/H)= \exp {\bf m}= \{ (\exp X, (\exp X)^{\lambda }); 
X \in sl_2(\mathbb R) \}$, 
\newline
where $\exp X \mapsto (\exp X)^{\lambda }: 
PSL_2(\mathbb R) \to PSL_2(\mathbb R)$ is a mapping. 
\end{Prop}
\Pro Let $S_X=\{ \exp t X; t \in \mathbb R \}$  be a $1$-parameter subgroup of  
$PSL_2(\mathbb R)$ isomorphic to $SO_2(\mathbb R)$. For all 
$x,y \in S_X$ is satisfied  
$(x y)^{\lambda}=x^{\lambda } y^{\lambda }$ and $(S_X ,S_X^{\lambda }) \cap H=\{(1,1) \}$. Hence the mapping $x \to x^{\lambda -1}$ is an automorphism of 
$S_X$. The only non-trivial automorphism of $S_X$ is the mapping 
$x \mapsto  x^{-1}$. Therefore the automorphism $x \mapsto x^{\lambda -1}$ 
must be the identity map and we have $\lambda =2$. 
For $x_1= \left ( \begin{array}{cc} 
 \frac{1}{2} & 0 \\
0 & 2  \end{array} \right )$ and $x_2= \left ( \begin{array}{rr} 
 \frac{1}{2} & -9 \\
0 & 2 \end{array} \right )$ we have 
\newline
\centerline{
$(x_i,x_i^2)=(R,1)(U_i,D^{-1} U_i D)$, } 
where 
$R=\left ( \begin{array}{cc} 
2 & 0 \\ 
0 &  \frac{1}{2}  \end{array} \right )$, $D= \left ( \begin{array}{cc} 
 \frac{\sqrt{5}}{5} & 0 \\ 
0 & \sqrt{5}  \end{array} \right )$,
and $U_i=D x_i^2 D^{-1}$. 
This means that the coset $(R,1) H^D$ of the conjugate subgroup $H^D$ of $H$ contains two different elements $(x_i, x_i^2)$ of 
$\sigma (G/H)$ $(i=1,2)$. Hence we have a contradiction to Proposition 1.6. in \cite{loops} (p. 19). 
\qed

\begin{Lemma} If the group  $G$ locally isomorphic to  $PSL_2(\mathbb C)$ is the group  topologically generated by the left 
translations of a $3$-dimensional almost differentiable left A-loop, then $G$ is isomorphic to $PSL_2(\mathbb C)$ and $H$ is 
isomorphic to $SO_3(\mathbb R)$. 
\end{Lemma}
\Pro According to  \cite{Asoh} (pp. 273-278) there are $4$ conjugacy classes of 
the $3$-dimensional subgroups of $G=SL_2(\mathbb C)$, which are denoted in \cite{Asoh} by 
$W_r$, $U_0$, $U_1$ and $SU_2(\mathbb C)$. 
Since the factor spaces $SL_2(\mathbb C)/U_i$ and $PSL_2(\mathbb C)/(U_i / \mathbb Z_2)$ for $i=0,1$ are homeomorphic to the 
topological direct product having as a factor the $2$-sphere or the projective plane respectively there is no differentiable loop 
realized on these factor spaces (cf. \cite{figula}, Proposition 2). 
\newline
Let now $H$ be locally isomorphic one of the subgroups $W_r$ or $W_r \mathbb Z_2/ \mathbb Z_2$, where 
\newline
\centerline{ $W_r=\left \{ \left(\begin{array}{cc}
\displaystyle \exp((ri-1)x) & 0 \\ z & \exp(-(ri-1)x)
\end{array} \right); x \in \mathbb R, z \in \mathbb C \right \}$ \quad for 
$r \in \mathbb R$. }
The Lie algebra ${\bf h}=w_r$ of the stabilizer $W_r$ has  following basis elements: 
\newline
\centerline{$ \{ r i K- K, i T- i U, U-T \}  \quad r \in \mathbb R $. }   
A  complement ${\bf m}$ to ${\bf h}$ in ${\bf g}$ contains a basis element  
$K+f (K)$ or $i K+ f(i K)$, where 
$f: {\bf m} \to {\bf h} $ is a linear map.  
Since the element 
\newline
\centerline{$ [U-T, K+f (K)]= [U-T, K] + [ U-T, f(K)] = $} 
\newline
\centerline{$ 2 U-2T +[U-T, f(K)] $} 
is an element of the intersection ${\bf h} \cap {\bf m}= \{ 0 \}$,  
 we have 
$[U-T, f(K)]=2 T-2 U$. This is the case precisely if  $f(K)=-K$
but then $f(K)$ is not an element of ${\bf h}$. This is a contradiction.  We obtain  the same contradiction if 
$i K+f(i K) \in {\bf m}$.  
\newline
Since $SU_2(\mathbb C)$ contains central elements $\neq 1$ of $SL_2(\mathbb C)$ the assertion follows.
\qed

\begin{Lemma} For all $a \in \mathbb R$ there is a reductive complement 
\newline
\centerline{
${\bf m}= \langle  T+a i T, i  U-a  U,  K+a i  K \rangle $ } 
to ${\bf h}=so_3(\mathbb R)$ generating ${\bf g}=sl_2(\mathbb C)$. 
\end{Lemma} 
\Pro  According to {\bf 1.2} let $\{ K, T, U, \hbox{i} K, \hbox{i} T, \hbox{i} U \}$ be a real basis of ${\bf g}=sl_2(\mathbb C)$. The Lie algebra ${\bf h}$ of the 
stabilizer $H=SO_3(\mathbb R)$ has the form ${\bf h}= \langle U, \hbox{i} T, \hbox{i} K \rangle $. An arbitrary component ${\bf m}$ to ${\bf h}$ has the shape
\newline
${\bf m}=\langle T+\hbox{a} U+\hbox{b i} T+\hbox{c i} 
K, 
\hbox{i} U+ \hbox{d} U+\hbox{e i} T+\hbox{f i} K,
K+ \hbox{g} U+\hbox{h i} T+ \hbox{k i} K \rangle $, where 
$ a,b,c,d,e,f,g,h,k \in \mathbb R$. The property $[{\bf h}, {\bf m}] \subseteq {\bf m}$ gives the assertion. 
\qed

\bigskip
\noindent
Now we determine the isomorphism classes and the isotopism classes of the loops $L_a$, $a \in \mathbb R$ belonging to the complements 
${\bf m}_a$.  
Two loops corresponding to $(G,H, \exp {\bf m}_a)$ and $(G,H, \exp {\bf m}_b)$ 
are isomorphic if and only if there exists an automorphism $\alpha $ of 
${\bf g}$ such that $\alpha ({\bf m}_a)={\bf m}_b$ and 
$\alpha ({\bf h})={\bf h}$. The automorphism group of ${\bf g}$ leaving 
${\bf m}_0$ and ${\bf h}$ invariant is the semidirect product $\Theta $ of 
$\hbox{Ad}_H$ and the group generated by the involutory map $\varphi :z \mapsto \bar{z}$. Since ${\bf m}$ is 
 a reductive subspace the  condition  $\alpha ({\bf m}_a)={\bf m}_b$, $\alpha \in \Theta $, is equivalent to   
$\varphi ({\bf m}_a)={\bf m}_b$. This identity is satisfied if and only if $b=-a$. Therefore a full isomorphism class consists of the loops $L_a$ and $L_{-a}$ $(a \in \mathbb R)$ and we may  choose as  representatives 
of these isomorphism classes the left A-loops $L_a,\ a \ge 0$. 
Since there is no $g \in G$ such that $g^{-1} {\bf m}_a g= {\bf m}_b$ for two different real numbers $a,b$ the isotopism classes and the isomorphism classes of the left A-loops $L_a, a \in \mathbb R$ are the same. 
\newline
The complement ${\bf m}_0=\langle T, \hbox{i} U, K \rangle $ satisfies $[{\bf m}_0, {\bf m}_0]={\bf h}$, 
 and 
$\bf{g}= {\bf m}_0 \oplus [{\bf m}_0,{\bf m}_0]$.
Hence it determines a $3$-dimensional connected Riemannian symmetric 
space (cf. \cite{ 
foundations1}, Chapter VI, Theorem 2.2 (iii)) and the loop $\hat{L}_0$ corresponding to the complement ${\bf m}_0$ is a Bruck loop. According to \cite{figula} this is the hyperbolic space loop.  
\newline
The loops $L_a$ for  $a \in \mathbb R$ have elementary models in the  upper half 
space $\mathbb R^{3+}=\{ (x,y,z) \in \mathbb R^3; z>0 \}$, which may be identified with the $\bf{J}$-quaternion 
space  (\cite{elementary}, p. 4) such that the point $j$ is the identity $e$ of $L_a$. The elements of  $L_a$ are the points of the 
$\bf{J}$-quaternion space. The  1-parameter subgroups through $e \in L_a$ have the same form for all 
loop 
$L_a$, but the sets $\exp {\bf m}_a$ differ. The multiplication in the loop $L_a$, 
$a \in \mathbb R$ is given by 
\newline
\centerline{$x \ast y= (\exp X) y$, for all $x,y \in L_a$}
where $X$ is the unique element of ${\bf m}_a$ such that $x = \exp  X$.  
\newline 
Summarizing our discussion we obtain

\begin{Theo}  If $L$ is a  connected almost differentiable strong\-ly left alternative left A-loop with dimension $3$ having a  semisimple 
Lie group $G$ as the group  
topologically generated by its left translations then $G$ is isomorphic to $PSL_2(\mathbb C)$ and the stabilizer $H$ of 
$e \in L$ in $G$ is isomorphic to $SO_3(\mathbb R)$ and $L=L_a$ is characterized by a real parameter $a$.  
\newline
The loops $L_a$ and $L_{-a}$ form a full isomorphism class, which is even a full isotopism class too. 
Among the loops $L_a$ only the hyperbolic space loop $L_0$ is a Bruck loop. This loop  is realized on the hyperbolic symmetric space by the multiplication 
$x \cdot y= \tau_{e,x} (y)$, where $\tau_{e,x}$ is the hyperbolic translation 
moving $e$ onto $x$. The tangent space ${\bf m}_0$ for the manifold of the left translations of $L_0$  is within 
the Lie algebra ${\bf g}$ of $G$ orthogonal to the  Lie algebra ${\bf h}$ 
of $H$ with respect to the Cartan-Killing form of ${\bf g}$.
\end{Theo}

\bigskip
\centerline{\bf 3-dimensional left A-loops corresponding to }
\centerline{\bf 4-dimensional non-solvable Lie groups}

\bigskip
\noindent
In this section we determine all 
 $3$-dimensional connected almost differentiable global left A-loops 
$L$ having a $4$-dimensional non-solvable Lie group $G$ as the group 
topologically 
generated by their left translations. Then the stabilizer $H$ of 
$e \in L$ 
in $G$ has dimension $1$.
\newline
In this case we have $G=PSL_2(\mathbb R) \times G_2$, where $G_2$ is one of 
the $1$-dimensional Lie groups, and $H$ is one of the cases 2 and 3 in the 
Proposition 5. 
\newline
The Lie algebra ${\bf g}$ of $G$ can be represented as $\bf{g}$$=sl_2(\mathbb R) \oplus  \mathbb R$. Let $(K,0)$, $(T,0)$, $(U,0)$ with $K,T,U$ defined in {\bf 1.2} be a real basis of $sl_2(\mathbb R) \oplus \{ 0 \}$ and let $(0,e_1)$ be the generator of $\{ 0 \} \oplus \mathbb R$. 

\begin{Lemma} The Lie algebra ${\bf g}$ is reductive with a $1$-dimensional subalgebra 
${\bf h}$ not contained in $sl_2(\mathbb R) \oplus \{ 0 \}$ and a $3$-dimensional complementary subspace ${\bf m}$ generating ${\bf g}$ in one of the following cases:
\newline
1) ${\bf h}=\langle (K, e_1) \rangle$, ${\bf m}_a=\langle (U,0), (T,0), (a K, (1+a) e_1) \rangle$, where $ a \in \mathbb R \backslash \{ -1 \} $
\newline
2) ${\bf h}=\langle (U+T, 2 e_1) \rangle$, ${\bf m}_b=\langle (U+T,0), (K,0), (U, 2 b e_1) \rangle$, where 
$ b \in \mathbb R \backslash \{ 0 \} $
\newline
3) ${\bf h}=\langle (U, e_1) \rangle$, ${\bf m}_c=\langle (K,0), (T,0), (c U, (1+c) e_1) \rangle$, where 
$ c \in \mathbb R \backslash \{ -1 \} $. 
\end{Lemma}
\Pro According to  Proposition 5 we may assume that ${\bf h}$ has one of the shapes given in $1$ till $3$. 
An arbitrary complement ${\bf m}$ to ${\bf h}$ in ${\bf g}$ has the shape in the case 1) 
\newline
\centerline{${\bf m}=\langle (U+a_1 K, a_1 e_1), (T+a_2 K, a_2 e_1), (a_3 K, (1+a_3) e_1) \rangle $}
in the case 2)  
\newline
${\bf m}=\langle (K + a_1(U+T),2 a_1 e_1), (U + a_2 (U+T),2 a_2 e_1), 
(a_3 (U+T), e_1 + 2 a_3 e_1) \rangle $
in the case 3) 
\newline
\centerline{${\bf m}=\langle (K + a_1\ U,a_1\  e_1), (T + a_2\  U, a_2\  e_1), 
(a_3\  U, e_1 + a_3\  e_1) \rangle ,$} 
where $a_1, a_2, a_3 \in \mathbb R$.
From the fact that $\hbox{dim}\  {\bf m} \cap (sl_2(\mathbb R) \oplus \{ 0 \})=2$ (Lemma 1) and from the property $[{\bf h}, {\bf m}] \subseteq {\bf m}$ follows the assertion. 
\qed 

\begin{Prop} The complement ${\bf m}_0$ in the case 3) is the unique reductive subspace such that there are global left A-loops $L_n$, $n \in \mathbb N$ 
with ${\bf m}_0=T_1 L_n$. These loops $L_n$ are Scheerer extensions of the Lie group $SO_2(\mathbb R)$ by the hyperbolic plane 
loop. 
\end{Prop}
\Pro  For all $a \in \mathbb R \backslash \{ -1 \}$ the complement ${\bf m}_a$ contains the elements 
\newline
\centerline{$k_a= -(1+d+d^2) (U,0)+ (-1+d+d^2) (T,0)+ \left ( \frac{a}{1+a} K,e_1 \right)$,  
for which } $Ad _g (k_a)=(K,e_1)$ holds with $(K,e_1) \in {\bf h}$ and  
$g=\left ( \pm \left (\begin{array}{cc}
1+d  & -1  \\
-d & 1 \end{array} \right), 0  \right) $ such that $d=- \frac{1}{2(1+a)}$.
This is a contradiction to Lemma 2. 
\newline
Now we deal with the complement ${\bf m}_b$. If $b > 0$ then for the elements 
\newline
$k_b=(K+U,2 b  e_1) \in {\bf m}_b$ are satisfied $Ad _g (k_b)=b (U+T, 2 e_1)$, where  
\newline
$b (U+T, 2 e_1) \in {\bf h}$ and  
$g=\left( \pm \left (\displaystyle \begin{array}{cc}
1 & - \frac{2 b}{\sqrt{2 b}} \\
\frac{1}{\sqrt{2 b}} &  0 \end{array} \right),0 \right) \in G$. This contradicts Lemma 2. 
\newline
For $b<0$  the subspace ${\bf m}_b$ contains the vectors 
\[ v_1= (-3 \pi b (U+T), 0), \quad v_2=(\sqrt{5 \pi ^2} K+ 3 \pi U, 6 \pi b e_1). \]
According to {\bf 1.2} the exponential images of the vectors $v_1$ and $v_2$ are  
\newline
\centerline{ $m_1=\exp {v_1}= \left (\left ( \begin{array}{cc}
1 & -6 \pi b \\
0 & 1 \end{array} \right ), 0 \right )$}   
and 
\newline
\centerline{ $m_2=\exp {v_2}= \left ( \pm I , 6 \pi b \right )$, } 
where $\pm I$ is the identity of $PSL_2(\mathbb R)$. 
One has $\left ( \pm I, 6 \pi b \right)=m_1 \cdot h_1=m_2$, where  $h_1= 
\left (\left ( \begin{array}{cc}
1 & 6 \pi b \\
0 & 1 \end{array} \right ), 6 \pi b \right )$. This is  a contradiction to Lemma 2. 
\newline
Finally we consider the reductive complements ${\bf m}_c$. 
For $c < -1$ and for the elements 
\newline
\centerline{
$k_c=\left ( \frac{1-2 e^4}{2 e^2} \right )(T,0)+ \left ( \frac{c}{1+c} U,e_1 \right ) \in {\bf m}_{c} $} 
we obtain $Ad_g (k_c)=(U,e_1)$, where $(U,e_1) \in {\bf h}$ and  
$g=\left (\pm \left (\begin{array}{cc}
 \frac{1}{e} & 0  \\
0 & e \end{array} \right), 0  \right)$ $\in G$, choosing $e$ such that 
$ \frac{1+2 e^4}{2 e^2} = \frac{c}{1+c}$. This is a contradiction to Lemma 2. 
\newline
If $c > -1$ but $c \neq 0$ the subspace ${\bf m}_c$ contains the vectors  
\newline
\centerline{$v_1=\left( k U,  \frac{k (1 +c)}{c} e_1 \right)$}
and 
\newline
\centerline{$v_2=\left( \sqrt{ \left ( \frac{ k^2 (1+c-nc)^2}{(1 +c)^2}-4 \pi ^2 \right) } T+  \frac{k (1+c-nc)}{1 +c} U, \frac{k (1+c-nc)}{c} e_1 \right)$.} 
According to 
{\bf 1.2} the images of $v_1, v_2$ under the exponential map have the forms: 
\newline
\centerline{$m_1=\exp {v_1}=
\left ( \pm \left ( \begin{array}{rr}
\cos {k} & \sin {k} \\
-\sin {k} & \cos {k} \end{array} \right ), 
 \left ( \begin{array}{rr}
\cos { \frac{k (1 +c)}{c}} & \sin { \frac{k (1 +c)}{c}} \\
-\sin { \frac{k (1 +c)}{c}} & \cos { \frac{k (1 +c)}{c}} \end{array} \right ) \right )$} 
and 
\newline
\centerline{$m_2=\exp {v_2}=\left ( \pm I , \left (\begin{array}{rr}
\cos { \frac{k (1+c-nc)}{c}} & \sin { \frac{k (1+c-nc)}{c} } \\
-\sin { \frac{k (1+c-nc)}{c} } & \cos { \frac{k (1+c-nc)}{c}} \end{array} \right ) \right )$ .} 
For  
\[ g= \left ( \pm I , \left (\begin{array}{rr}
\cos { \frac{k (1+c-nc)}{c}} & \sin { \frac{k (1+c-nc)}{c} } \\
-\sin { \frac{k (1+c-nc)}{c} } & \cos { \frac{k (1+c-nc)}{c}} \end{array} \right ) \right ) \in G ,\] where $k \in \mathbb Z $ such 
that $k > 
\sqrt {\frac{4 \pi ^2(1+c)^2}{(1+c-nc)^2}}$ and $\pm I$ is the identity of $PSL_2(\mathbb R)$, one has $g=m_1 \cdot h_1=m_2$ such that 
\newline
\centerline{$h_1=\left ( \pm \left ( \begin{array}{rr}
\cos {k} & -\sin {k} \\
\sin {k} & \cos {k} \end{array} \right ), 
 \left ( \begin{array}{rr}
\cos {n k} & -\sin {n k} \\
\sin {n k} & \cos {n k} \end{array} \right ) \right )$. }  
This again contradicts Lemma 2. 
\newline
For $c=0$ the complement ${\bf m}_{c}$ has the shape:
${\bf m}_0 = \langle (K,0), (T,0), (0, e_1) \rangle $. 
\newline
Since $\big [ [{\bf m}_0,{\bf m}_0], {\bf m}_0 \big ] \subseteq {\bf m}_0$ the loops $L$ with the property $T_1 L={\bf m}_0$ are 
global  Bol loops.  According to $\cite{figula}$ we have a global Bol loop $L_n$ for all $n \in \mathbb N$ having the direct product 
$PSL_2(\mathbb R) \times SO_2(\mathbb R)$ as the group topologically generated by its left translations and as the stabilizer $H$ of 
$e \in L_n$ in $G$ the group $H_n=\{ ( x,x^n) | x \in SO_2(\mathbb R), n \in \mathbb N \}$. The non-isotopic loops $L_n$ are 
 Scheerer extensions of the Lie group $SO_2(\mathbb R)$ by the hyperbolic plane loop (cf. \cite{loops}, Section 2).  \qed

\bigskip
\noindent
The loops $L_n$, $n \in \mathbb N$ are homeomorphic to $G/H_n$, which is the cylinder $\mathbb R^2 \times S^1$. Let $\tilde L$ be the 
universal covering of $L$. Since  $\tilde L$ is homeomorphic to $\mathbb R^3$ the loop   $\tilde L$ contains a central subgroup 
isomorphic to $\mathbb Z$. Moreover all other coverings of $L$ is $\tilde L / n \mathbb Z$. The universal covering group $\tilde G$ 
of $G$ is the group 
$\widetilde{PSL_2(\mathbb R)} \times \mathbb R$, which contains the central subgroup $\pi_1(G)=\mathbb Z \times \mathbb Z$. The 
universal 
covering group $\tilde H$ of $H$ is the group $\{ (x, n x) | x \in \mathbb R, n \in \mathbb N \}$, which is isomorphic to 
$\mathbb R$. The group $G^{\ast }$ topologically generated by the left translations of $\tilde{L}$ is the covering group 
$\tilde{G}/ \{ (z, nz) | z \in \mathbb Z, n \in \mathbb N \}$ and the stabilizer of the identity of $\tilde{L}$ is the group 
$\tilde{H} \pi_1(G)/ \pi_1(G)$.     
\newline
Summarizing our discussion we have: 

\begin{Theo}
There are precisely three  classes ${\cal C}_1$, 
${\cal C}_2$ and ${\cal C}_3$  of  connected strongly left alternative almost differentiable left A-loops with dimension $3$ such that the group $G$ 
topologically generated by their  left translations is a $4$-dimensional non-solvable Lie group.  
\newline
The class ${\cal C}_1$ consists of left A-loops $L$ such that the group $G$ topologically generated by their left translations 
is  isomorphic to $PSL_2(\mathbb R) \times \mathbb R$ and  the stabilizer of the identity of these loops is isomorphic to 
$SO_2(\mathbb R) \times \{ 0 \}$.   
Every loop in  ${\cal C}_1$ is a 
product of a $2$-dimensio\-nal loop isomorphic to the hyperbolic plane loop 
 with the Lie group $\mathbb R $. These loops are not isotopic. The only differentiable Bruck loop in ${\cal C}_1$  corresponds 
to the 
section $\sigma :G/H \to G$ such that $\sigma (G/H)=M_1 \times SO_2(\mathbb R)$, where $M_1$ is the image of the section of the 
hyperbolic plane.   
\newline
In the class ${\cal C}_2$ are the products of a $2$-dimensio\-nal loop isomorphic to the hyperbolic plane loop 
 with the Lie group $SO_2(\mathbb R)$.  These loops are not isotopic and  the group $G$ topologically generated by their left 
translations is  isomorphic to 
$PSL_2(\mathbb R) \times SO_2(\mathbb R)$ and the stabilizer of the identity of these loops is isomorphic to 
$SO_2(\mathbb R) \times \{ 1 \}$. In ${\cal C}_2$ there is again precisely one differentiable Bruck loop $\hat L$. The 
image of the 
section of $\hat L$ is the direct product of the image of the hyperbolic plane loop with the Lie group $SO_2(\mathbb R)$.  
\newline
In the class ${\cal C}_3$ are contained  the Scheerer extensions $L_n$, $n \in \mathbb N $ of the Lie group 
$SO_2(\mathbb R)$ by the 
hyperbolic plane loop and the coverings of $L_n$.   The group $G$ topologically generated by the left translations of $L_n$ is the direct 
product $PSL_2(\mathbb R) \times SO_2(\mathbb R)$ and the stabilizer $H$ of $e \in L_n$ in $G$ is the group 
$H_n=\{ (x, x^n), x \in SO_2(\mathbb R), n \in \mathbb N \}$. 
\newline
The intersection of the classes ${\cal C}_2$ and  ${\cal C}_3$ is the Bruck loop $\hat L$.  
\end{Theo}

\bigskip
\centerline{\bf 3-dimensional left A-loops belonging to } 
\centerline{\bf 5-dimensional non-solvable Lie groups}

\bigskip
\noindent
Now we determine the $3$-dimensional connected almost differentiable global 
left A-loops
having a $5$-dimensional non-solvable Lie group $G$ as the group 
topologically generated by the left translations of $L$. 
In this case the stabilizer of $e \in L$ in $G$ is a $2$-dimensional
closed subgroup of $G$ containing no non-trivial normal subgroup of $G$.
Then because of Proposition 5 we have to 
investigate only the following case:
\newline
$G$ is locally isomorphic to the semi-direct product $PSL_2(\mathbb R) \ltimes \mathbb R^2$, which 
is the connected component of the group for area preserving affinities of 
$\mathbb R^2$. 

\bigskip
\noindent
For the Lie algebra ${\bf g}=sl_2(\mathbb R) \ltimes \mathbb R^2$ of $G$ we can choose the following basis 
elements: 
\newline
\centerline{$K=\left( \begin{array}{rrr}
0 & 0 & 0 \\
0 & 1 & 0 \\
0 & 0 & -1 \end{array} \right)$,\ $T=\left( \begin{array}{ccc}
0 & 0 & 0 \\
0 & 0 & 1 \\
0 & 1 & 0 \end{array} \right)$,\ $U=\left( \begin{array}{rrr}
0 & 0 & 0 \\
0 & 0 & 1 \\
0 & -1 & 0 \end{array} \right)$, }
\newline
\centerline{$ e_1=\left( \begin{array}{ccc}
0 & 1 & 0 \\
0 & 0 & 0 \\
0 & 0 & 0 \end{array} \right)$, \  $e_2=\left( \begin{array}{ccc}
0 & 0 & 1 \\
0 & 0 & 0 \\
0 & 0 & 0 \end{array} \right)$. }
The multiplication table is given by:
\newline
\centerline{$[K,e_1]=[T,e_2]=-[U, e_2]=-2\ e_1,\ [K,e_2]=-[U,e_1]=-[T,e_1]=2 \ e_2, $}
\newline
\centerline{$[e_1,e_2]=0, \ [K,T]=2\ U,\  [K,U]=2\ T,\ [U,T]=2\ K. $}

\begin{Lemma} The Lie algebra ${\bf g}=sl_2(\mathbb R) \ltimes \mathbb R^2$ is reductive with a subalgebra ${\bf h}$ which does not contain any ideal $\neq 0$ of ${\bf g}$ and a $3$-dimensional complementary subspace ${\bf m}$ generating ${\bf g}$ in the following case  
${\bf h}=\langle K, e_1 \rangle $ and  
\newline
${\bf m}= \langle e_2, U+b e_1, T- b e_1 \rangle $, where $b \in \mathbb R$. 
\end{Lemma}
\Pro  The $2$-dimensional Lie algebras ${\bf h}$ of ${\bf g}$, which does not contain any ideal $\neq 0$ of ${\bf g}$ (up to mapping 
$Ad_g, g \in G$) are 
${\bf h}_1=\langle K, U-T \rangle $, ${\bf h}_2=\langle K, e_1 \rangle $, ${\bf h}_3=\langle U-T, e_1 \rangle $. We have for a complement ${\bf m}$ to $\bf{h}_1$ in $\bf{g}$ the general form: 
\newline
${\bf m}=\langle e_1+a_1 K + a_2 (U-T), e_2 + b_1 K +b_2 (U-T), 
U+ c_1 K+ c_2 (U-T) \rangle $.  
A complement ${\bf m}$ to $\bf{h}_2$ in $\bf{g}$ we can write in the following form:
\newline
\centerline{${\bf m}=\langle e_2+a_1 K+a_2 e_1, U +b_1 K +b_2 e_1, T + c_1 K +c_2 e_1 
\rangle $.} 
An arbitrary complement $\bf{m}$ to $\bf{h}_3$ in $\bf{g}$ can be given as 
follows: 
\newline
${\bf m}= \langle e_2+a_1(U-T)+a_2 e_1,\ K+b_1(U-T)+b_2 e_1,\ U+c_1(U-T)+
c_2 e_1 \rangle $. 
Here $a_1,a_2,b_1,b_2,c_1,c_2$ are real parameters.  The assertion follows now  from the property 
$[{\bf h}, {\bf m} ] \subseteq {\bf m} $. \qed

\begin{Prop} There is no global left A-loop $L$ corresponding to the reductive subspace 
${\bf m}= \langle e_2, U+b e_1, T- b e_1 \rangle $, $b \in \mathbb R$.
\end{Prop} 
\Pro The element $e_2 \in {\bf m}$ is equal to $Ad_g (e_1)$, where  $e_1  \in {\bf h}$ and  
$g= \left ( \begin{array}{rrr}
1 & 0 & 0 \\
0 & 0 & 1 \\
0 & -1 & 0 \end{array} \right ) \in G$. This is a contradiction to Lemma 2. 
\qed 
\newline
This consideration yields the following
\begin{Theo} There is no $3$-dimensional connected almost differentiable 
global left A-loop $L$ having a $5$-dimensional non-solvable Lie group as the 
group topologically generated by its left translations.
\end{Theo}

\newpage
\centerline{\bf 3-dimensional left A-loops with } 
\centerline{\bf 6-dimensional non-solvable Lie groups}

\bigskip
\noindent 
Now we  determine all $3$-dimensional connected almost differentiable left A-loops such that the group $G$ 
topologically generated by their left translations is a non-semisimple and non-solvable  Lie group.  
According to  Lemma 4 and  
Propositions 5 we have to discuss the following cases 
\newline
$\alpha $) $G$ is locally isomorphic to $PSL_2(\mathbb R) \ltimes \mathbb R^3$,
\newline
$\beta $) $G$ is the group for orientation preserving affinities of $\mathbb R^2$, 
\newline
$\gamma $) $G$ is locally isomorphic to $SO_3(\mathbb R) \ltimes \mathbb R^3$, which is the connected component of the euclidean motion group of $\mathbb R^3$.  

\bigskip
\noindent
In the case $\alpha $) the group multiplication in $G$ is 
given by 
\[ (A_1, X_1) \circ (A_2, X_2)= (A_1\ A_2,\ A_2^{-1}\ X_1\ A_2+ X_2), \]
where $(A_i,X_i)$, $i=1,2$ are  two elements of $G$ such that $X_i$ $(i=1,2)$ 
are represented by $(2 \times 2)$ real  matrices with trace $0$.
\newline
A basis of the Lie algebra ${\bf g}=sl_2(\mathbb R) \ltimes \mathbb R^3$ of $G$ can be chosen as follows:
\newline
\centerline{$e_1=\left( 0, \left( \begin{array}{rr}
0 & -1 \\
1 & 0 \end{array} \right) \right)$, $e_2=\left( \left( \begin{array}{rr}
1 & 0 \\
0 & -1 \end{array} \right), 0 \right)$, $e_3=\left( \left( \begin{array}{cc}
0 & 1 \\
1 & 0 \end{array} \right), 0 \right)$,}   
\newline
\centerline{$ e_4=\left( \left( 
\begin{array}{rr}
0 & 1 \\
-1 & 0 \end{array} \right), 0 \right)$, 
$e_5=\left( 0, \left( \begin{array}{rr}
-1 & 0 \\
0 & 1 \end{array} \right) \right)$, $e_6=\left( 0, \left( 
\begin{array}{cc}
0 & 1 \\
1 & 0 \end{array} \right) \right)$. } 
According to \cite{jacobson} (p. 17) we obtain the following 
 multiplication table in ${\bf g}$: 
\newline
\centerline{$[e_1,e_2]=:e_6,\ [e_1,e_3]=:e_5,\ [e_2,e_3]=:e_4,\ [e_5,e_4]=-e_6, $}
\newline
\centerline{$[e_1,e_4]= [e_1,e_5]=[e_1,e_6]= [e_2,e_5]= [e_3, e_6]=[e_6,e_5]=0, $} 
\newline
\centerline{$[e_2,e_6]= [e_3,e_5]= -e_1, \ [e_2,e_4]=e_3,\ [e_3,e_4]=-e_2,\  
[e_6,e_4]=e_5. $} 

\begin{Lemma} The Lie algebra ${\bf g}=sl_2(\mathbb R) \ltimes \mathbb R^3$ is reductive with a subalgebra ${\bf h}$ containing no non-zero ideal of ${\bf g}$ and a $3$-dimensional complementary subspace ${\bf m}$ generating ${\bf g}$ in one of the following cases:  
\newline
(i) ${\bf h}=\langle e_1, e_2, e_6 \rangle $ and  ${\bf m}_{b_2,b_3}= \langle e_5, e_3-b_3 e_1-b_2 e_6, e_4+ b_2 e_1+b_3 e_6 \rangle $, where $b_2, b_3 \in \mathbb R$. 
\newline
(ii) ${\bf h}=\langle e_2, e_3, e_4 \rangle $ and  ${\bf m}_{a}= \langle e_1+a e_4, e_6-a  e_3, e_5+ a e_2 \rangle $, where $a \in \mathbb R \backslash \{ 0 \}$. 
\newline 
(iii) ${\bf h}=\langle e_4,e_5,e_6 \rangle $ and ${\bf m}_{b_1,b_2}= \langle e_1, e_2+b_1 e_6+b_2 e_5, e_3-b_2 e_6+b_1 e_5 \rangle $, where $b_1,b_2 \in \mathbb R$.
\end{Lemma}
\Pro The $3$-dimensional subalgebras ${\bf h}$  of ${\bf g}$, which does not contain any non zero-ideal are the following:
\newline
a) $\langle e_2,\ e_5,\ e_1+e_6 \rangle $,
\newline
b) $ \langle e_2+k\ e_5,\ e_1,\ e_6 \rangle$, where $k \in \mathbb R$,
\newline
c) $\langle e_3+e_4,\ e_5,\ e_1-e_6 \rangle $,
\newline 
d) $\langle e_2,\ e_3+e_4,\ e_1-e_6 \rangle $,    
\newline
e) $\langle e_2,e_3,e_4 \rangle $,
\newline
f) $\langle e_4,e_5,e_6 \rangle $.     
\newline
In the case a) the basis elements of an arbitrary complement ${\bf m}$ to 
${\bf h}$ in ${\bf g}$ are: 
\newline
\centerline{$f_1=e_1+a_1 e_2+a_2 e_5+a_3 (e_1+e_6)$ , $f_2=e_4+ b_1 e_2+b_2 e_5+b_3 (e_1+e_6)$,  }
\newline
\centerline{$f_3=e_3+c_1 e_2+c_2 e_5+c_3 (e_1+e_6), $}
with $a_1,a_2,a_3,b_1,b_2,b_3,c_1,c_2,c_3 \in \mathbb R$.
\newline
If ${\bf h}=\langle e_2+k\ e_5,\ e_1,\ e_6 \rangle$ with $k \in \mathbb R$
then an arbitrary complement ${\bf m}$ to ${\bf h}$ has as basis elements: 
\newline
\centerline{$f_1= e_5+ a_1 (e_2+k e_5)+ a_2 e_1 +a_3 e_6$, $f_2= e_3 +c_1 (e_2+k e_5)+ c_2 e_1+ c_3 e_6$,    }
\newline
\centerline{$f_3= e_4+ b_1( e_2+k e_5)+b_2 e_1 +b_3 e_6 ,$ }
where $a_1,a_2,a_3,b_1,b_2,b_3 ,c_1,c_2,c_3$ are real parameters. 
\newline
In the case c) we can choose as basis elements of an arbitrary complement ${\bf m}$ to 
${\bf h}$  the following: 
\newline
\centerline{$f_1=e_1 +a_1 (e_3+e_4)+ a_2 e_5 +a_3 (e_1-e_6), $} 
\newline
\centerline{$f_2= e_2 +b_1 (e_3+e_4)+b_2 e_5 +b_3 (e_1-e_6), $}
\newline
\centerline{$f_3=e_3+ c_1 (e_3+e_4) +c_2 e_5 +c_3 (e_1-e_6), $} 
where $a_1,a_2,a_3,b_1,b_2,b_3,c_1,c_2,c_3 \in \mathbb R$. 
\newline
In the case d) the generators of an arbitrary complement ${\bf m}$ to 
${\bf h}$ in ${\bf g}$ are: 
\newline
\centerline{$f_1= e_1+a_1 e_2+ a_2( e_3+e_4) + a_3(e_1-e_6), $}
\newline
\centerline{$f_2= e_5+b_1 e_2+ b_2( e_3+e_4) + b_3( e_1-e_6), $}
\newline
\centerline{$f_3= e_3+ c_1 e_2+ c_2(e_3+e_4)+ c_3(e_1-e_6), $}
with the real parameters $a_1,a_2,a_3,b_1,b_2,b_3,c_1,c_2,c_3$.  
\newline
In the case e) one has ${\bf h} \cong sl_2(\mathbb R)$. An arbitrary complement 
${\bf m}$ to ${\bf h}$ has the form 
\newline
\centerline{$\langle e_1 + a_1 e_2 + a_2 e_3 + a_3 e_4, e_6 + b_1 e_2 + b_2 e_3 + b_3 e_4, e_5 + c_1 e_2 + c_2 e_3 + c_3 e_4 \rangle$, }  
with $a_1,a_2,a_3,b_1,b_2,b_3,c_1,c_2,c_3 \in \mathbb R$. 
\newline
Now we consider the last case. An arbitrary complement ${\bf m}$ to 
${\bf h}$ in ${\bf g}$ has the following basis elements: 
\newline
\centerline{$ \{ e_1+a_1 e_4+a_2 e_5+a_3 e_6, e_2+b_1 e_4+b_2 e_5+b_3 e_6, 
e_3+c_1 e_4+c_2 e_5+c_3 e_6 \} , $} where  
$a_1,a_2,a_3,b_1,b_2,b_3,c_1,c_2,c_3 \in \mathbb R$. 
Using the relation $[ {\bf h}, {\bf m}] \subseteq {\bf m}$ we obtain the assertion. 
\qed

\begin{Prop} There is no global left A-loop $L$ belonging to the reductive subspaces (i) and (ii) in Lemma 18.  
\end{Prop}
\Pro The element $e_5 \in {\bf m}$ in the case (i) is equal to  $Ad_g (e_6)$, such that $e_6 \in {\bf h}$ and  
$g=\left ( \pm \left( \begin{array}{rr}
 \frac{1}{2} &  - \frac{1}{2} \\
1 & 1 \end{array} \right), 0 \right) \in G$.  
The element $e_6-a e_3+e_1+a e_4 \in {\bf m}_a$ in the case (ii) 
 for all $a \in \mathbb R \backslash \{ 0 \}$ is equal to  $Ad_g (a (e_4-e_3))$ with $a (e_4-e_3) \in {\bf h}$ 
and $g= \left (1, \left( \begin{array}{cc}
 \frac{1}{2 a} & 0 \\
0 & - \frac{1}{2 a} \end{array} \right) \right) \in G$. These facts contradict Lemma 2. \qed

\bigskip
\noindent
Now we deal with the case (iii) in Lemma 18.
Since the group $SL_2(\mathbb R)$ has no $3$-dimensional linear representation the group $G$ is isomorphic to the semidirect 
product of $PSL_2( \mathbb R) \ltimes \mathbb R^3$  
and $H$ is  
isomorphic to the following $3$-dimensional subgroup of $G$: 
\newline
\centerline{$H=\left \{ \left( \pm \left( \begin{array}{rr}
\cos {t} & \sin {t} \\
- \sin {t} & \cos {t} \end{array} \right),  \left( \begin{array}{rr}
-x & y \\
y  & x  \end{array} \right) \right), t \in [0, 2 \pi ), x,y \in \mathbb R 
\right \}$. } 
\newline
\newline
Now we determine the isomorphism classes and the isotopism classes of the left A-loops $L_{b_1,b_2}$ having the subspaces ${\bf m}_{b_1,b_2}$ $(b_1,b_2 \in \mathbb R)$ as the tangent spaces $T_1 L_{b_1,b_2}$. 
\newline
We have precisely two isomorphism classes ${\cal C}_i$ $(i=1,2)$ of the loops 
$L_{b_1,b_2}$  belonging to the triples $(G,H,\exp {\bf m}_{b_1,b_2})$ for all $b_1,b_2 \in \mathbb R$. 
\newline
The first class ${\cal C}_1$ consists loops belonging to ${\bf m}_{b_1,b_2}$ for $b_2=0$. 
Denote by $\hat{{\bf m}}_{b_1}$ the complement ${\bf m}_{b_1,0}$ for all $b_1 \in \mathbb R$.
One has 
$[\hat{{\bf m}}_{b_1}, \hat{{\bf m}}_{b_1}]={\bf h}$ and  
${\bf g}=\hat{{\bf m}}_{b_1} \oplus 
[\hat{{\bf m}}_{b_1}, \hat{{\bf m}}_{b_1}]$ for all $b_1 \in \mathbb R$. 
Every loop $L_{b_1,0}$ in ${\cal C}_1$ is a Bruck loop and as a representative of this class we may choose the loop 
$\hat{L_0}=L_{0,0}$. According to \cite{figula} the loop $\hat{L_0}$ is a global differentiable Bruck loop, which is called the pseudo-euclidean space loop. 
\newline
The other class ${\cal C}_2$ consists of loops $L_{b_1,b_2}$ having 
$T_1 L_{b_1,b_2}={\bf m}_{b_1,b_2}$ for $b_2 \neq 0$. 
Since the automorphism $\beta $ of the Lie algebra ${\bf g}$ defined by 
\newline
$ \beta \  (e_1)=\sqrt {c^2+d^2}\  e_1$, 
\newline
$\beta \  (e_6)=-d\ e_5 + c\  e_6$,
\newline
$\beta \  (e_5)= c \ e_5 + d\ e_6$,
\newline
$\beta \ (e_4)=e_4 $,
\newline
$\beta \ (e_2)=\displaystyle \frac{c}{\sqrt {c^2+d^2}}\ e_2 
- \frac{d}{\sqrt {c^2+d^2}}\ e_3 
-c\ b_1\ e_6 +d\ b_1\ e_5 $,
\newline
$\beta \ (e_3)=\displaystyle  \frac{c}{\sqrt {c^2+d^2}}\ e_3
+\ \frac{d}{\sqrt {c^2+d^2}}\ e_2- d \ b_1\  e_6\ - c\ b_1\ e_5 $, 
\newline
where $\varepsilon \sqrt {c^2+d^2}= \displaystyle \frac {1}{b_2}$ with $\varepsilon =1$ for $b_2 >0$  and $\varepsilon =-1$ for $b_2 <0$,  
 leaves the subalgebra ${\bf h}$ invariant and   
$\beta ({\bf m}_{b_1,b_2})={\bf m}_{0,1}$ for all $b_1 \in \mathbb R, b_2 \in \mathbb R \backslash \{ 0 \}$ holds, we may  
choose the loop ${\hat L}_1=L_{0,1}$ as a representative of the class ${\cal C}_2$. 
\newline
Since there is no $g \in G$ such that $g^{-1} {\bf m}_{b_1,b_2} g = {\bf m}_{b_1',0}$ holds with $b_1, b_1' \in \mathbb R$, $b_2 \in \mathbb R \backslash \{ 0 \}$ 
the isotopism classes of the left A-loops $L_{b_1,b_2}$ coincide with the isomorphism classes ${\cal C}_1$, ${\cal C}_2$. 
\newline
Now we prove that  ${\hat L}_1=L_{0,1}$ is a global left A-loop. 
\newline
The exponential map $\exp : {\bf g} \to G$ is described in section 7 in \cite{figula}.  
\newline
The image of ${\bf m}_{0,1}$ under the exponential map is given as follows:  
The subspace ${\bf m}_{0,1}$ has the shape 
\newline
\centerline{$ {\bf m}_{0,1}= \left \{ \left ( \left ( \begin{array}{rr}
\lambda _2 & \lambda_3  \\
 \lambda _3 & -\lambda _2 \end{array} \right), \left( \begin{array}{cc}
-\lambda _2 & -\lambda _1 -\lambda_3 \\
\lambda _1 -\lambda _3 & \lambda _2
\end{array} \right) \right); \lambda_1, \lambda_2, \lambda_3 \in \mathbb R \right \}$.} 
According to {\bf 1.2} the first component of $ \exp {\bf m}_{0,1}$
is 
\newline
\centerline{$\left (\pm \left( \begin{array}{cc}
\cosh \sqrt {A}+  \frac{ \sinh \sqrt {A}}{\sqrt {A}} \lambda_2 & 
 \frac{ \sinh \sqrt {A}}{\sqrt {A}} \lambda _3 \\
 \frac{ \sinh \sqrt {A}}{\sqrt {A}}\lambda _3  
& \cosh \sqrt {A}- \frac{ \sinh \sqrt {A}}{\sqrt {A}} \lambda_2 \end{array} 
\right) \right)$, } 
the second component of  $\exp {\bf m}_{0,1}$  is
$\left( \begin{array}{rr}
r'(1) & s'(1) \\
v'(1) & -r'(1) \end{array} \right), $ where 
\[ r'(1)=\frac{\lambda _3\  \lambda _1}{4 A}\ 
(e^{ \sqrt {A}}- e^{ - \sqrt 
{A}})^2-\lambda_2 , \] \[ s'(1)= \frac{- \lambda _1}
{4 \sqrt {A} } (e^{2 \sqrt {A}}- e^{ -2 \sqrt 
{A}}) -\frac{\lambda _2 \lambda _1}{4 A}
(e^{ \sqrt {A}}- e^{ - \sqrt {A}})^2 - \lambda _3, \]
\[ v'(1)= \frac{\lambda _1}
{4 \sqrt {A} }(e^{2 \sqrt {A}}- e^{ -2 \sqrt {A}}) -\frac{\lambda _2 \lambda _1}{4 A}
(e^{ \sqrt {A}}- e^{ - \sqrt {A}})^2 -\lambda _3,  \] and 
$A=\lambda _2^2+\lambda _3^2$. 
\newline 
The submanifold $\exp {\bf m}_{0,1}$ is the image of a sharply transitive global section $\sigma :G/H \to G$ if and only if each element $g \in G$ can be uniquely written as a product $g=m h$ with $m \in \exp {\bf m}_{0,1}$ and $h \in H$, moreover $ \exp {\bf m}_{0,1}$ operates sharply transitively on $G/H$. 
\newline
In \cite{figula} section 7 we have shown that each element of $G=PSL_2(\mathbb R)$ can  be uniquely written  as
\newline
\centerline{
$\left( \pm \left( \begin{array}{cc}
a & b \\
c & d  \end{array} \right), \left( \begin{array}{rr}
x & y \\
z & -x \end{array} \right) \right)= $ }
\newline
\centerline{$ \left( \left( \begin{array}{ll}
a_1 & 0 \\
b_1 & a_1^{-1} \end{array} \right),  \left( \begin{array}{rr}
0 & u \\
-u & 0 \end{array} \right) \right) \cdot  \left( \pm \left( \begin{array}{rr}
\cos {t} & \sin {t} \\
- \sin {t} & \cos {t} \end{array} \right), \left( \begin{array}{rr}
k & l \\
l & -k \end{array} \right) \right) $}
with $a,b,c,d \in \mathbb R, \ a d-b c=1,\ x,y,z \in \mathbb R, a_1>0, b_1 \in \mathbb R,\ t \in [0,2 \pi)$,
such that $k=x, l=\displaystyle \frac{y+z}{2}, u=\displaystyle \frac{y-z}{2}$. 
Therefore it is sufficient to prove that there is 
to each element $g \in G$ with the shape 
\[  \left(  \left( \begin{array}{ll}
a & 0 \\
b & a^{-1} \end{array} \right),  \left( \begin{array}{rr}
0 & u \\
-u & 0 \end{array} \right) \right); a>0,b,u \in \mathbb R  \] 
precisely one  
$m \in \exp {\bf m}_{0,1}$ and $h \in H$ such that $g=m\ h$ or  equivalently $m=g\ h^{-1}$. 
\newline
The first component of $\exp {\bf m}_{0,1}$ is precisely  
the section $\sigma _1$ of the hyperbolic plane loop given in \cite{loops} (pp. 281-282).  Therefore  
for given $a>0,b \in \mathbb R$ we have unique 
$ \lambda_2, \lambda_3 \in \mathbb R, t \in [0, 2 \pi)$ such that 
\newline
\centerline{$\left( \begin{array}{cc}
 \cosh \sqrt {A}+  \frac{ \sinh \sqrt {A}}{\sqrt {A}} \lambda_2 & 
 \frac{ \sinh \sqrt {A}}{\sqrt {A}} \lambda _3 \\
\frac{ \sinh \sqrt {A}}{\sqrt {A}}\lambda _3  
& \cosh \sqrt {A}- \frac{ \sinh \sqrt {A}}{\sqrt {A}} 
\lambda_2 \end{array} 
\right)= $}
\newline
\centerline{$\left( \begin{array}{ll}
a & 0 \\
b & a^{-1} \end{array} \right) \left( \begin{array}{rr}
\cos {t} & \sin {t} \\
- \sin {t} & \cos {t} \end{array} \right)$, }  
where $A=\lambda _2^2+\lambda _3^2$. Hence we have to consider the second component 
 of $\exp {\bf m}_{0,1}$.  For given $u, \lambda_2, \lambda_3$ we 
have to find unique $\lambda_1,k,l \in \mathbb R$ such that 
\[ \left( \begin{array}{rr}
r'(1) & s'(1) \\
v'(1) & -r'(1) \end{array} \right)= \left( \begin{array}{cc}
k & l+u \\
l-u & -k \end{array} \right), \]
where 
$r'(1), s'(1), v'(1)$ are values of functions, which depend on 
the variables 
$\lambda_1, \lambda_2,\lambda_3$. 
Since for $\lambda_1$ we obtain the equation 
\newline
\centerline{$2u= \displaystyle \frac{-\lambda_1}{2 \sqrt{\lambda_2^2+\lambda_3^2}} 
(e^{2 \sqrt{ \lambda_2^2+\lambda_3^2}}-e^{-2 \sqrt{\lambda_2^2+\lambda_3^2}})$}
we have for the unique solutions  
$\lambda_1=\displaystyle \frac{-4 u \sqrt {\lambda_2^2+\lambda_3^2}}
{\displaystyle 
e^{2 \sqrt{ \lambda_2^2+\lambda_3^2}}-e^{-2 \sqrt{\lambda_2^2+\lambda_3^2}}}$, $k=r'(1)$,  
$l=\displaystyle \frac {s'(1)+v'(1)}{2}$. 
\newline
Now we verify that the section $\sigma _1$ corresponding to the loop ${\hat L}_1$ is sharply transitive, this means that for given elements  
\[ \left(  \left( \begin{array}{ll}
a_1 & 0 \\
b_1 & a_1^{-1} \end{array} \right),  \left( \begin{array}{cc}
0 & u_1 \\
-u_1 & 0 \end{array} \right) \right) \ \ \hbox{and} \ \  \left(  \left( \begin{array}{ll}
a_2 & 0 \\
b_2 & a_2^{-1} \end{array} \right),  \left( \begin{array}{cc}
0 & u_2 \\
-u_2 & 0 \end{array} \right) \right), \]
where  $a_1>0,a_2>0,b_1,b_2,u_1,u_2 \in \mathbb R$ there exists precisely one element $z \in \exp {\bf m}_{0,1}$ 
and a  $h= \left( \pm \left( \begin{array}{rr}
\cos {t} & \sin {t} \\
- \sin {t} & \cos {t} \end{array} \right), \left( \begin{array}{rr}
k & l \\
l & -k \end{array} \right) \right) \in H$, where $t,k,l \in \mathbb R$ such that the equation 
\newline
(I) \centerline{ $z \left(  \left( \begin{array}{ll}
a_1 & 0 \\
b_1 & a_1^{-1} \end{array} \right),  \left( \begin{array}{rr}
0 & u_1 \\
-u_1 & 0 \end{array} \right) \right)= $}
\newline
\centerline{$\left(  \left( \begin{array}{ll}
a_2 & 0 \\
b_2 & a_2^{-1} \end{array} \right),  \left( \begin{array}{rr}
0 & u_2 \\
-u_2 & 0 \end{array} \right) \right) \left( \pm \left( \begin{array}{rr}
\cos {t} & \sin {t} \\
- \sin {t} & \cos {t} \end{array} \right), \left( \begin{array}{rr}
k & l \\
l & -k \end{array} \right) \right) $}
holds.  
The real variables $\lambda_1, \lambda_2, \lambda_3$ of $z \in \exp {\bf m}_{0,1}$ are determined by the following equations
\newline
1. \centerline{$\displaystyle \frac{ \sinh \sqrt{A}}{\sqrt{A}} \left( \lambda_2 \left( a_1+\displaystyle \frac{a_2^2}{a_1} \right) + \lambda_3 \left( b_1+\displaystyle \frac{b_2 a_2}{a_1} \right) \right)+ $}
\newline
\centerline{$\cosh \sqrt{A} \left ( a_1-\displaystyle \frac{a_2^2}{a_1} \right) =0 $}
\newline
2. \centerline{$ \displaystyle \frac{ \sinh \sqrt{A}}{\sqrt{A}} \left( \lambda_2 \left ( 
\displaystyle \frac{b_2 a_2}{a_1} -b_1 \right) + \lambda_3 \left ( 
\displaystyle \frac{a_1^2+b_2^2}{a_1} \right) \right)+ $} 
\newline
\centerline{$\cosh \sqrt{A} \left (b_1 - \displaystyle \frac{b_2 a_2}{a_1} \right)=0 $} 
3. \centerline{$ 2(u_2-u_1)+\lambda_3( b_1^2-a_1^2+a_1^{-2})+2 a_1 b_1 \lambda_2= $}
\newline
\centerline{$\displaystyle \frac{- \lambda_1 (b_1^2+a_1^2-a_1^{-2})}{4 A} (e^{2 \sqrt {A}}- e^{ -2 \sqrt 
{A}})+  $}
\newline
\centerline{$\lambda_1 
\left( \displaystyle  \frac{\lambda_2(a_1^2-b_1^2-a_1^{-2})+ 2 \lambda_3 b_1 a_1}{4 A} \right)
(e^{ \sqrt {A}}- e^{ - \sqrt {A}})^2, $} 
where $A=\lambda _2^2+\lambda _3^2$.
\newline
If $z$ is an element of ${\bf m}_{0,0}$ in the equation (I)
then we obtain for the variables $\lambda_1, \lambda_2, \lambda_3$ of $z \in \exp {\bf m}_{0,0}$ the above equations 1, 2, and the equation 
\newline
3'. \centerline{$2(u_2-u_1)=
 \displaystyle \frac{- \lambda_1 (b_1^2+a_1^2-a_1^{-2})}{4 A} (e^{2 \sqrt {A}}- e^{ -2 \sqrt 
{A}})+ $}
\newline
\centerline{$ \lambda _1 \left( \displaystyle  \frac{\lambda_2(a_1^2-b_1^2-a_1^{-2})+ 2 \lambda_3 b_1 a_1}{4 A} \right)
(e^{ \sqrt {A}}- e^{ - \sqrt {A}})^2$ .}  
The equations 1, 2, 3', have unique solutions because $\sigma _0$ is a sharply transitive section. Therefore the equations 1, 2, 3, are also uniquely solvable for the variables $\lambda_1, \lambda_2, \lambda_3$. 
Hence the sharply transitive global section $\sigma _1$ yields also a global loop ${\hat L}_1(\sigma _1)$, which is a proper left 
A-loop.

\bigskip
\noindent 
An elementary model of the loop ${\hat L}_1$ may be  given on the set $\Psi $ of the euclidean planes in the pseudo-euclidean affine 
space (cf. \cite{vorlesungen}).  
The elements of the loops ${\hat L}_1$ are the same as the elements of  ${\hat L}_0$ (\cite{figula}), but the sets of the left 
translations 
 $\exp {\bf m}_{0,1}$ respectively  $\exp {\bf m}_{0,0}$ and hence the multiplication of these two loops differ.  The 
multiplication  in the loop ${\hat L}_1$ is given by 
\newline
$(\ast \ast)$ \centerline{$ \quad \quad Q_1 \ast Q_2= \tau _{P,Q_1} (Q_2)$, \quad for all 
$Q_1, Q_2 \in \Psi $,  } 
\newline
where $\tau_{P,Q_1}$ is the unique element of $\exp {\bf m}_{0,1}$ mapping 
the plane $P$, which is the identity of ${\hat L}_1$ onto $Q_1$.

\bigskip
\noindent
In the case $\beta $) the Lie algebra ${\bf g}$ of the group $G$  has a  real basis  
\newline
\centerline{$ e_1=\left( \begin{array}{ccc}
0 & 0 & 0 \\
0 & 1 & 0 \\
0 & 0 & 0 \end{array} \right)$,\ $e_2=\left( \begin{array}{ccc}
0 & 0 & 0 \\
0 & 0 & 1 \\
0 & 0 & 0 \end{array} \right)$,\ $e_3=\left( \begin{array}{ccc}
0 & 0 & 0 \\
0 & 0 & 0 \\
0 & 1 & 0 \end{array} \right)$, }
\newline
\centerline{$e_4=\left( \begin{array}{ccc}
0 & 0 & 0 \\
0 & 0 & 0 \\
0 & 0 & 1 \end{array} \right)$,\ $e_5=\left( \begin{array}{ccc}
0 & 1 & 0 \\
0 & 0 & 0 \\
0 & 0 & 0 \end{array} \right)$,\ $e_6=\left( \begin{array}{ccc}
0 & 0 & 1 \\
0 & 0 & 0 \\
0 & 0 & 0 \end{array} \right)$. }
The multiplication is given by the following rules:
\newline
\centerline{$[e_1,e_2]=[e_2,e_4]=e_2,\ [e_1,e_3]=[e_3,e_4]=-e_3,\ [e_1,e_5]=[e_3,e_6]=-e_5, $} 
\newline
\centerline{$[e_1,e_6]=[e_1,e_4]=[e_2,e_6]=[e_3,e_5]=[e_4,e_5]=[e_5,e_6]=0, $}
\newline
\centerline{$[e_2,e_3]=e_1-e_4,\ [e_2,e_5]=[e_4,e_6]=-e_6. $}

\begin{Lemma} The Lie algebra ${\bf g}$ is reductive with a subalgebra ${\bf h}$ containing no non-zero ideal of ${\bf g}$ and a $3$-dimensional complementary subspace ${\bf m}$ generating ${\bf g}$ in the following case:  
 ${\bf h}=\langle e_1, e_4, e_5 \rangle $ and  ${\bf m}= \langle e_2, e_3, e_6 \rangle $. 
\end{Lemma}
\Pro The $3$-dimensional subalgebras ${\bf h}$, which does not contain any ideal $\neq 0$ of ${\bf g}$ are 
\newline
a) ${\bf h}=\langle e_1-e_4, e_2, e_3 \rangle $
\newline
b) ${\bf h}=\langle e_1, e_2, e_4 \rangle $
\newline
c) ${\bf h}=\langle e_1, e_4, e_5 \rangle $
\newline
d) ${\bf h}=\langle e_1+e_4, e_3, e_5 \rangle $
\newline
e) ${\bf h}=\langle e_3, e_4, e_5 \rangle $
\newline
f) ${\bf h}=\langle e_1-e_4, e_3, e_5 \rangle $. 
\newline
The basis elements of an arbitrary complement 
${\bf m}$ to ${\bf h}$ in ${\bf g}$ in the case a) are: 
\newline
\centerline{$ e_1+a_1(e_1-e_4)+a_2 e_2+a_3 e_3$,  $e_5+b_1(e_1-e_4) +b_2 e_2+ b_3 e_3$,    } 
\newline
\centerline{$ e_6+c_1(e_1-e_4)+c_2 e_2+c_3 e_3, $}
where $a_1,a_2,a_3,b_1,b_2,b_3,c_1,c_2,c_3 $ are real parameters.  
\newline
In the case b) an arbitrary complement ${\bf m}$ to ${\bf h}$ in ${\bf g}$ 
has the following shape: 
\newline
\centerline{
${\bf m}= \langle e_3+a_1 e_1+a_2 e_2+a_3 e_4,\  e_5 +b_1 e_1 +b_2 
e_2+ b_3 e_4,$ }
\newline 
\centerline{$e_6+c_1 e_1+c_2 e_2+c_3 e_4 \rangle $, } 
with the real numbers 
$a_1,a_2,a_3,b_1,b_2,b_3,c_1,c_2,c_3 \in \mathbb R$. 
\newline
In the case c) an arbitrary complement ${\bf m}$ to ${\bf h}$ in ${\bf g}$ 
has as generators: 
\newline
\centerline{$\{ f_1=e_2+a_1 e_1+a_2 e_4+a_3 e_5,\  f_2=e_3 +b_1 e_1 +b_2 
e_4+ b_3 e_5,$} 
\newline
\centerline{$f_3=e_6+c_1 e_1+c_2 e_4+c_3 e_5 \}, $} where 
$a_1,a_2,a_3,b_1,b_2,b_3,c_1,c_2,c_3 \in \mathbb R$. 
\newline
In the case d) the basis elements of an arbitrary complement ${\bf m}$ to ${\bf h}$ in ${\bf g}$ are: 
\newline
\centerline{$ e_1+a_1(e_1+e_4)+a_2 e_3+a_3 e_5$,  $e_2 +b_1(e_1+e_4) +b_2 e_3+ b_3 e_5$,  }
\newline
\centerline{$  e_6+c_1(e_1+e_4)+c_2 e_3+c_3 e_5 , $} with the real parameters
$a_1,a_2,a_3,b_1,b_2,b_3,c_1,c_2,c_3$.  
\newline
In the case e) an arbitrary complement ${\bf m}$ to ${\bf h}$ in 
${\bf g}$ is given by: 
\newline
\centerline{$\langle e_1+a_1 e_3+a_2 e_4+a_3 e_5,\  e_2 +b_1 e_3 +b_2 
e_4+ b_3 e_5,\  e_6+c_1 e_3+c_2 e_4+c_3 e_5 \rangle , $} with 
$a_1,a_2,a_3,b_1,b_2,b_3,c_1,c_2,c_3 \in \mathbb R$. 
\newline
In the last case f) for the basis elements of an arbitrary complement 
${\bf m}$ to ${\bf h}$ in ${\bf g}$ one has: 
\newline
\centerline{$ e_1+a_1(e_1-e_4)+a_2 e_3+a_3 e_5$,\   $e_2 +b_1(e_1-e_4) +b_2 e_3+ b_3 e_5$,    }
\newline
\centerline{$  e_6+c_1(e_1-e_4)+c_2 e_3+c_3 e_5 , $}
where $a_1,a_2,a_3,b_1,b_2,b_3,c_1,c_2,c_3$ are real numbers.
The assertion follows from the property $[{\bf h}, {\bf m}] \subseteq {\bf m}$. 
\qed 

\begin{Prop}
There is no global left A-loop $L$ having the reductive subspace ${\bf m}= \langle e_2, e_3, e_6 \rangle $ as the tangent space $T_1 L$. 
\end{Prop}
\Pro The subspace ${\bf m}$ contains the element 
$e_2+e_3$, which is equal to   $Ad _g (e_1-e_4)$, where $e_1-e_4  \in {\bf h}$ and 
 $g= \left (\displaystyle \begin{array}{rrr} 
1 & 0 & 0 \\
0 & 1 &  \frac{1}{2} \\
0 & -1 & \frac{1}{2} \end{array} \right) \in G$. 
This is a contradiction to Lemma 2. 
\qed

\bigskip
\noindent
Now we consider the case that $G$ is locally isomorphic to $SO_3(\mathbb R) 
\ltimes \mathbb R^3$. This group can be represented by the pairs of complex 
$(2 \times 2)$-matrices 
\newline
\centerline{
$(A,X)=\left ( \pm \left ( \begin{array}{rr}
a & b \\
-\bar{b} & \bar{a} \end{array} \right ), \left( \begin{array}{cc}
k & l i+n \\
-l i + n & -k \end{array} \right ) \right); $} 
$a,b \in \mathbb C, a \bar{a}+b \bar{b}=1, k,l,n \in \mathbb R$.  Here $\bar{a}$ denotes  the complex conjugate of $a \in \mathbb C$. 
The group multiplication is 
given by the rule 
\newline
\centerline{$ (A_1, X_1) \circ (A_2, X_2)= (A_1\ A_2,\ A_2^{-1}\ X_1\ A_2+ X_2)$. }
Any $3$-dimensional subgroup $H$ of $G=SO_3(\mathbb R) \ltimes \mathbb R^3$,  
which contains no  non-trivial normal subgroup of $G$, 
is locally isomorphic either to a semidirect product of a $2$-dimensional translation group by a 
$1$-dimensional rotation group  $SO_2(\mathbb R)$ or to the subgroup $\{(a,0), a \in SO_3(\mathbb R) \}$. The 
first possibility cannot occur since the factor space $G/H$ is the topological product having as a factor the 
$2$-sphere which is not parallelizable.       
\newline
Now we deal with the second possibility for $H$.

\begin{Lemma} For all $a \in \mathbb R \backslash \{ 0 \}$ there is a reductive complement 
\newline
\centerline{${\bf m}_a=\langle V_1 +a Z, V_2 + a Y, V_3 -a X \rangle $}
 to the Lie algebra ${\bf h}_2$ of 
$H_2$ generating ${\bf g}=so_3(\mathbb R) \ltimes \mathbb R^3$. 
\end{Lemma}
\Pro 
Denote by $X,Y,Z$ the generators correspond to $1$-dimensional rotations and 
let $V_3,V_2,V_1$ be the axes of the rotation groups corresponding to $X,Y$ 
respectively $Z$. We can identify the basis elements of ${\bf g}$ with the following matrices:
\newline
\centerline{$ X= \left( \left( \begin{array}{rr}
i & 0  \\
0 & -i \end{array} \right),0 \right)$,\ $Y= \left( \left( \begin{array}{cc}
0 & i  \\
i & 0  \end{array} \right),0 \right)$, \ $Z= \left( \left( \begin{array}{rr}
0 & -1  \\
1 & 0 \end{array} \right),0 \right)$, } 
\newline
\centerline{$V_1=  \left( 0, \left( \begin{array}{rr}
0 & i \\
-i & 0 \end{array} \right) \right)$,\ $V_2= \left( 0, \left( \begin{array}{cc}
0 & 1  \\
1 & 0  \end{array} \right) \right)$,\ $V_3= \left( 0, \left( \begin{array}{rr}
-1 & 0  \\
0 & 1  \end{array} \right) \right)$. }
According to \cite{jacobson} (p. 17) 
the multiplication table of ${\bf g}= su_2(\mathbb C) \ltimes \mathbb R^3$ 
is given by:
\newline
\centerline{$[X,Y]=Z,\ [Z,X]=Y,\ [Y,Z]=X,\ [X,V_1]=[Z,V_3]=-V_2, $} 
\newline
\centerline{$[X,V_2]=[Y,V_3]=V_1,\  [Z,V_2]=-[Y,V_1]=V_3, $}
\newline
\centerline{$[X,V_3]=[Y,V_2]=[V_1,V_2]=[V_1,V_3]=[V_2,V_3]=[Z,V_1]=0. $}
The Lie algebra ${\bf h}_2$ of $H_2$ has as generators $X, Y, Z$.
An arbitrary complement ${\bf m}$ to ${\bf h}_2$ in ${\bf g}$ has the following 
shape: 
\newline
\centerline{${\bf m}= \langle V_1+a X+b Y+c Z,\ V_2+d X+e Y+f Z,\ V_3+g X+h Y+i Z 
\rangle ,$ }
where $a,b,c,d,e,f,g,h,i \in \mathbb R$. The subspace ${\bf m}$ 
satisfies the condition 
\newline
$[{\bf h}_2, {\bf m}] \subseteq {\bf m}$ if and only if 
${\bf m}$ has the following form:  
\newline
\centerline{${\bf m}_a= \langle V_1 +a Z, V_2 + a Y, V_3 -a X \rangle , $}
where  $a \in \mathbb R \backslash \{ 0 \}$. 
\qed

\bigskip
\noindent
Using the automorphism $\varphi $ of ${\bf g}$ having the form: 
\newline
$\varphi (V_1)=  \frac{1}{2 a}\  V_1 +  \frac{\sqrt{3} }{2 a} \  V_2$, 
\newline
$\varphi (V_2)=  \frac{\sqrt{3} }{2 a}\  V_1 - \frac{1}{2 a}\  V_2$,  
\newline
$\varphi (V_3)=-  \frac{1}{ a}\ V_3$, 
\newline
$\varphi (X)=- X$, 
\newline
$\varphi (Y)=-  \frac{1}{2}\  Y +  \frac{\sqrt{3}}{2}\  Z$,
\newline
$\varphi (Z)=  \frac{\sqrt{3}}{2}\ Y +   \frac{1}{2}\  Z$, 
\newline
 for all 
$a \in \mathbb R \backslash \{ 0 \}$ we have 
$\varphi ({\bf h})= {\bf h}$  and $\varphi ({\bf m}_a)= {\bf m}_1$. Therefore the local loops $L_a$ having $T_1 L_a={\bf m}_a$ form an isomorphism class $\cal{C}$ and as a representative 
of $\cal{C}$ can be chosen the local loop $L_1$ belonging to  ${\bf m}_1$.

\begin{Prop} The local loop $L_1$ is not a global left A-loop. 
\end{Prop}
\Pro The exponential map $\exp : {\bf g} \to G$ is given by the following way: For
$X \in {\bf g}$ we have $\exp X= v _X(1)$, where $v _X(t)$ is the 
1-parameter subgroup of $G$ with the property 
$ \frac{d}{dt} \big |_{t=0} v _X(t)=X$. 
In the 1-parameter subgroup $ \alpha (t)=( \beta (t), \gamma (t))$ of $G$
with the conditions 
\newline
\centerline{
$\alpha (t=0)=(1,0)$ and   
$ \frac{d}{dt} \big |_{t=0} \alpha (t)=(X_1,X_2)=X \in \bf{g}$ } 
the first component 
$\beta (t)$ is the 1-parameter subgroup of 
$SO_3( \mathbb R)$ and the second component $\gamma (t)$ satisfies 
\[ \displaystyle \frac{d}{dt} \gamma (t)= \frac{d}{ds} \big| _{s=0} \gamma (t+s)=
- \frac{d}{ds} \big| _{s=0} \beta (s) \gamma(t)+\gamma(t)  
\frac{d}{ds} \big| _{s=0} \beta (s)+ \frac{d}{ds} \big| _{s=0} \gamma (s)= \]
\[ -X_1 \gamma(t) + \gamma(t) X_1 +X_2. \] 
For $X_1=\left (  \begin{array}{cc}
\lambda _1 i  & \lambda _2 i - \lambda _3 \\
-\lambda _2 i +\lambda _3 & -\lambda _1 i \end{array} \right )$, $ X_2= \left (  \begin{array}{cc}
\lambda _5  & \lambda _4 i + \lambda _6 \\
-\lambda _4 i +\lambda _6 & -\lambda _5 \end{array} \right )$ and $\gamma (t)= \left( \begin{array}{cc}
  r(t)  & v(t) i + s(t) \\
 - v(t) i + s(t)  &  - r(t)  \end{array} \right)$, where  $\lambda _1,\lambda _2, \lambda _3, \lambda _4,\lambda _5, \lambda _6 \in \mathbb R$  we get the following inhomogen system of linear differential equations:  
\[  \frac{d}{dt} \left ( \begin{array}{c}
r(t) \\
s(t) \\
v(t) \end{array} \right )= \left ( \begin{array}{ccr}
0 & -2 \lambda _1 & -2 \lambda _3 \\
2 \lambda _1 & 0 & 2 \lambda _2   \\
2 \lambda _3 & -2 \lambda _2 & 0 \end{array} \right ) \left ( \begin{array}{c}
r(t) \\
s(t) \\
v(t) \end{array} \right ) + \left ( \begin{array}{c}
\lambda _5 \\
\lambda _6 \\
\lambda _4 \end{array} \right ) \]
with the following initial conditions: 
\[ r(0)=s(0)=v(0)=0,\  \frac{d}{dt} \big| _{t=0} r(t)=\lambda _5 , \ 
\frac{d}{dt} \big| _{t=0} s(t)=\lambda _6,\ \frac{d}{dt} \big| _{t=0} v(t)=\lambda _4.\] 
The solution of this inhomogeneous 
system is:
\[ r(t)= - \frac{i[(e^{2 lit}-e^{2 lit})(\lambda _5 \lambda _1^2+ \lambda _5 \lambda _3^2-\lambda _6 \lambda _3 \lambda _2+\lambda _4 \lambda _1 \lambda _2)]}{4 l^3} \]
\[-\frac{\big [ \big ( e^{ lit}-e^{- lit}
\big )^2 (-\lambda _6 \lambda _1 - \lambda _4 \lambda _3)+t(4 \lambda _4 \lambda _1 \lambda _2-4 \lambda _6 \lambda _3 \lambda _2-4 \lambda _5 \lambda _2^2) \big ] }{4 l^2} , \]
\[ s(t)= - \frac{i (e^{2 lit}-e^{-2 lit})(\lambda _6 \lambda _1^2+ \lambda _6 \lambda _2^2-\lambda _5 \lambda _3 \lambda _2+\lambda _4 \lambda _1 \lambda _3) }{4 l^3} \]
\[-\frac{\big [ \big ( e^{ lit}-e^{-lit}
\big )^2 (\lambda _4 \lambda _2 + \lambda _5 \lambda _1)+t(4 \lambda _4 \lambda _1 \lambda _3-4 \lambda _5 \lambda _3 \lambda _2-4 \lambda _6 \lambda _3^2) \big ]}{4 l^2} , \]
\[ v(t)= - \frac{i (e^{2 lit}-e^{-2 lit})(\lambda _4 \lambda _3^2+ \lambda _4 \lambda _2^2+\lambda _5 \lambda _1 \lambda _2+\lambda _6 \lambda _1 \lambda _3)}{4 l^3} \]
\[-\frac{\big [ \big ( e^{ lit}-e^{- lit}
\big )^2 (\lambda _5 \lambda _3 - \lambda _6 \lambda _2)+t(4 \lambda _5 \lambda _1 \lambda _2+4 \lambda _6 \lambda _3 \lambda _1-4 \lambda _4 \lambda _1^2) \big ] }{4 l^2} , \] where $l=\sqrt {\lambda _1^2+\lambda _2^2+\lambda _3^2}$. 
Since ${\bf m}_1$ has the form 
\[ {\bf m}_1= \left \{ \left( \left ( \begin{array}{cc}
-c i & -a+b i \\
a+b i & c i \end{array} \right ) , \left ( \begin{array}{cc} 
-c  & a i+b  \\
-a i+b  & c  \end{array}  \right ) \right ); a,b,c \in \mathbb R \right \}, \] 
according to {\bf 1.2} (Section 1) 
the first component of $ \exp {\bf m}_1$ is 
\[ (\exp {\bf m}_1)_1= \left ( \begin{array}{cc}
\cos{\sqrt{k}}- \frac{c i \sin{\sqrt{k}}}{\sqrt{k}} &  \frac{(-a+b i) \sin{\sqrt{k}}}{\sqrt{k}} \\ 
 \frac{(a+b i) \sin{\sqrt{k}}}{\sqrt{k}} & \cos{\sqrt{k}}+ \frac{c i \sin{\sqrt{k}}}{\sqrt{k}} \end{array} \right ), \] 
the second component of $\exp {\bf m}_1$ has the shape   
\newline
\centerline{
$(\exp {\bf m}_1)_2=\left( \begin{array}{cc}
r(1) & v(1)i+s(1) \\
-v(1)i+s(1) & -r(1) \end{array} \right)$, }  
where    
\newline
\centerline{$r(1)=-c ( e^{ \sqrt{k}i}-e^{- \sqrt{k}i})^2$,  $s(1)=b ( e^{ \sqrt{k}i}-e^{- \sqrt{k}i})^2 $, 
$v(1)=a ( e^{ \sqrt{k}i}-e^{- \sqrt{k}i})^2 $, } 
and  $k=a^2+b^2+c^2$.  
\newline
From the equation $ g=\left (1, \left ( \begin{array}{cc}
0 & f i \\
-f i & 0 \end{array} \right ) \right )=(( \exp {\bf m}_1)_1, ( \exp {\bf m}_1)_2) (h,0)$  
with $f \neq 0$ one has $h=( \exp {\bf m}_1)_1^{-1}$. This means that  
\newline
\centerline{$\left ( \begin{array}{cc}
0 & f i  \\
-f i & 0 \end{array} \right ) = $} 
\newline
\centerline{$(\exp {\bf m}_1)_1 \left (\begin{array}{cc} 
-(e^{ k i}-e^{- k i})^2 c & (e^{ k i}-e^{- k i})^2 (a i +b) \\
(e^{ k i}-e^{- k i})^2 (-a i +b) & (e^{ k i}-e^{- k i})^2 c \end{array} \right ) (\exp {\bf m}_1)_1^{-1}, $} where $k=\sqrt{a^2+b^2+c^2}$. Hence we obtain 
\newline
\centerline{$-c ( e^{ \sqrt{a^2+b^2+c^2}i}-e^{- \sqrt{a^2+b^2+c^2}i})^2=0$,\   $ a ( e^{ \sqrt{a^2+b^2+c^2}i}-e^{- \sqrt{a^2+b^2+c^2}i})^2=f$, }
\newline
\centerline{$ b ( e^{ \sqrt{a^2+b^2+c^2}i}-e^{- \sqrt{a^2+b^2+c^2}i})^2=0$ . }
\newline
Therefore we may assume that $c=b=0$. Then we have 
\newline
\centerline{$a ( e^{ \sqrt{a^2}i}-e^{- \sqrt{a^2}i})^2 = f \ \hbox{or} \ a(\sinh{\sqrt{a^2}i})^2=-a (\sin{\sqrt{a^2}})^2=\displaystyle \frac{f}{4}. $} 
Since the function $x \mapsto  -x (\sin{\sqrt{x^2}})^2$ is not injective, there exist  different real numbers $a_1,a_2$ with the properties $\sin{(\sqrt{ a_1^2})} \neq \sin{(\sqrt{ a_2^2})}$ but $a_1(\sin{\sqrt{a_1^2}})^2= a_2 (\sin{\sqrt{a_2^2}})^2$.  Hence 
$g$ can be written in different way as a product of an element in $\exp {\bf m}_1$ and an element of $H$ which contradicts    
Lemma 2.   
\qed

\noindent
From the above discussion we obtain: 
\begin{Theo}
There is only one class ${\cal C}$ of the $3$-dimensional connected almost differentiable left A-loops $L$ such that the group $G$ topologically generated by the left translations  is a $6$-dimensional non semisimple and non-solvable Lie group. 
The group  $G$ is isomorphic to the semidirect product $PSL_2(\mathbb R) \ltimes \mathbb R^3$,  where 
the action of $PSL_2(\mathbb R)$ on $\mathbb R^3$ is the adjoint action of 
$PSL_2(\mathbb R)$ on its Lie algebra, and the stabilizer of the identity of the loops in ${\cal C}$ is the $3$-dimensional subgroup of $G$ 
\newline
\centerline{$\left \{ \left(\pm \left( \begin{array}{rr}
\cos t & \sin t \\
-\sin t & \cos t \end{array} \right) , \left( \begin{array}{rr}
-x & y \\
y & x \end{array} \right) \right); t \in [0 , 2 \pi ),x,y \in \mathbb R 
\right \}. $}
The loops in the class ${\cal C}$ can be characterized by two real parameters $a,b$ and form precisely two isomorphism classes 
which coincide the isotopism classes. The one isomorphism class containing  the Bruck loops 
$L_{b_1,0}$, $b_1 \in \mathbb R$ has as a representative  the pseudo-euclidean space loop $L_{0,0}={\hat L}_0$. 
As a representative of the other isomorphism class consisting  of left A-loops $L_{b_1,b_2}$ with $b_2 \neq 0$ may be chosen the 
 loop $L_{0,1}={\hat L}_1$. The loops ${\hat L}_0$ and  ${\hat L}_1$ are realized on the pseudo-euclidean affine 
space $E(2,1)$.  The elements of these loops are the planes on which the euclidean metric is induced but the sets of left translations 
differ.  
\end{Theo}

\bigskip
\noindent
$\begin{array}{lcl}
\hbox{{\small \'Agota Figula}} & \quad & \hbox{{\small and}} \\
\hbox{{\small Mathematisches Institut }} & \quad &  \hbox{{\small  Institute of Mathematics }} \\
\hbox{{\small der Universit\"at Erlangen-N\"urnberg}} & \quad & \hbox{{\small University of Debrecen}} \\
\hbox{{\small  Bismarckstr. 1 $\frac{1}{2}$, }} & \quad & \hbox{{\small  P.O.B. 12, H-4010 Debrecen, }} \\ 
\hbox{{\small D-91054 Erlangen, Germany,}} & \quad & \hbox{{\small Hungary,}} \\
\hbox{{\small figula@mi.uni-erlangen.de}} & \quad & \hbox{{\small  figula@math.klte.hu}} \end{array}$

\end{document}